\documentclass[twoside]{article}
\usepackage[a4paper]{geometry}
\usepackage[latin1]{inputenc} % ou \usepackage[utf8]{inputenc}
\usepackage[T1]{fontenc} % ou \usepackage[OT1]{fontenc}
\usepackage{RR}
\usepackage{hyperref}
\usepackage{amsfonts}
\usepackage{amsmath}
\usepackage{amssymb}
\usepackage{pdfpages}
\usepackage{amsthm,euscript,makeidx,pifont,boxedminipage,subeqnarray}
\usepackage{subfig}
\usepackage{pdfpages}
\usepackage{amsmath,amssymb,amsfonts,amsthm,euscript,makeidx,pifont,boxedminipage}
\usepackage{color}
\usepackage{empheq}
\usepackage{epsfig}
\usepackage{graphicx}

\newcommand{\R}{\mathbb{R}}
\newcommand{\w}{\Omega}
\newcommand{\woned}{\Omega_{1D}}
\newcommand{\wtwod}{\Omega_{2D}}
\newcommand{\wtwo}{\Omega_{2}}
\newcommand{\frad}{\displaystyle\frac}
\newcommand{\dn}[1]{ \frad{\partial {#1}}{\partial n }}
\newcommand{\dx}[1]{ \frad{\partial {#1}}{\partial x }}	
	
\newcommand{\dxx}[1]{ \frad{\partial^2 {#1}}{\partial x^2 }}	

\newcommand{\dnone}[1]{ \frad{\partial {#1}}{\partial n_1 }}
\newcommand{\dntwo}[1]{ \frad{\partial {#1}}{\partial n_2 }}
	
\newcommand{\be}{\begin{equation}}
\newcommand{\ee}{\end{equation}}
\newcommand{\av}[1]{ \overline{#1}}

\newtheorem{theorem}{Theorem}

\newtheorem{remark}{Remark}
\newtheorem{definition}{Definition}

\newtheorem{prop}{Proposition}
\graphicspath{{Figs/}}
\RRNo{8182}
\RRdate{December 2012}
\RRauthor{% les auteurs
M. Tayachi\thanks{Inria, Laboratoire Jean Kuntzmann, BP 53, 38041 Grenoble Cedex 9, France.} , 
A. Rousseau \thanks{Inria, Laboratoire Jean Kuntzmann, BP 53, 38041 Grenoble Cedex 9, France.} , 
E. Blayo\thanks{Universit\'e de Grenoble \& Inria, Laboratoire Jean Kuntzmann, BP 53, 38041 Grenoble Cedex 9, France.} , 
N. Goutal\thanks{EDF R\&D and Laboratoire d'hydraulique Saint-Venant.} and 
V. Martin\thanks{LAMFA UMR-CNRS 7352, Universit\'e de Picardie Jules Verne, 33 Rue St. Leu, 80039 Amiens, France.}
}
\RRtitle{Une m\'ethode de couplage  de type Schwarz dans le cadre d'un probl\`eme multi-dimensionnel}
\RRetitle{Design and analysis of a Schwarz coupling method for a dimensionally heterogeneous problem}
\titlehead{DDM for a dim-heterogeneous problem}
\RRresume{Dans ce document nous \'etudions et analysons et une m\'ethode de couplage multi-dimensionnel it\'erative. Nous consid\'erons 
le cas de l'\'equation de Laplace 2-D avec des conditions aux bords non sym\'etriques, coupl\'ee avec une \'equation de Laplace 1-D correspondante. 
dans un premier temps nous montrons comment obtenir le mod\`ele 1-D \'a partir du mod\`ele 2-D par int\'egration verticale et par analogie avec la 
d\'erivation des \'equations de Saint-Venant \'a partir des \'equations de Navier-Stokes. Ensuite nous pr\'esentons un algorithme de couplage de type Schwarz.
Nous discutons le choix des conditions aux interfaces de couplage. Nous d\'emontrons la convergence de tels algorithmes et donnons quelques r\'esultats 
th\'eoriques sur le choix de la position des interfaces de couplage. Un r\'esultat th\'eorique sur le contr\^ole de l'erreur entre la solution globale 2-D de r\'ef\'erence  et la solution 2-D coupl\'ee sera aussi donn\'e. Enfin nous illustrons ces r\'esultats num\'eriquement.
}
\RRabstract{
 In the present work, we study and analyze an efficient  iterative coupling method for a dimensionally heterogeneous problem . We
consider the case of 2-D Laplace equation with non symmetric boundary conditions with
 a corresponding 1-D Laplace equation.  We will first show how to obtain the 1-D model from the 2-D one by integration along
 one direction, by analogy with the link between shallow water equations and the Navier-Stokes system.
  Then, we will focus on the design of an 
 Schwarz-like iterative coupling method. We will discuss  the choice of boundary conditions at coupling interfaces. 
 We will prove the convergence of such algorithms and   give some theoretical  results related to the choice of the location of the coupling interface, and the control
 of the difference between a global 2-D reference solution and the 2-D coupled one. These theoretical results will be illustrated numerically. }
\RRmotcle{couplage multi-dimensionnel, d\'ecomposition de domaine, analyse muli-\'echelles}
\RRkeyword{dimensionally heterogeneous coupling, domain decomposition methods, multiscale analysis}
\RRprojets{MOISE}
 \RCGrenoble % Grenoble - Rh\^one-Alpes
\begin{document}
\newcommand{\saut}{\vspace*{5mm}\par}
\newcommand{\id}{I\hspace*{-0.6mm}d}
\makeRR   % cas d'un rapport de recherche
%% \makeRT % cas d'un rapport technique.
%% a partir d'ici, chacun fait comme il le souhaite

\section{Introduction}
 Hydrodynamical phenomena can be described by a wide variety of mathematical and numerical models, spanning a large range of possible levels of complexity and realism. When dealing with the representation of a complex fluid system, such as an ensemble of rivers and channels or a human blood system, the dynamical behavior of the flow is often spatially heterogeneous. This means that it is generally not necessary to use the most complex model everywhere, but that one can adapt the choice of the model to the local dynamics. One has then to couple several different models, corresponding to different areas. Such an approach is generally efficient from a computational point of view, since it avoids heavy computations with a full complex model in areas where a simpler model is able to represent the dynamics quite accurately. Thus this makes it possible to build a hybrid numerical representation of an entire complex system, while its simulation with a unique model would be either non relevant with a simple one or too expensive with a complex one.\par
 In such a hierarchy of models, the simplest ones are often simplifications of the more complex ones. Let mention for instance the so called ``primitive equations", which are widely used to represent the large scale ocean circulation, and are obtained by making some assumptions in the Navier-Stokes equations. It is important to note that such simplifications may involve a change in the geometry and in the dimension of the physical domain, thus leading to simplified models which are $m$-D while the original one was $n$-D, with $n>m$.
An obvious example is given by the shallow water equations, which are derived from the Navier-Stokes equations by integration along the vertical axis, see for instance 
\cite{gerbeauperthame} for a rigorous mathematical derivation using asymptotic analysis techniques in the 2-D to 1-D case
%, or  \cite{UNKNOWN} \color{red} [\`a trouver] \color{black} for a derivation in a more physical perspective)
. Such a coupling between dimensionally heterogeneous models has been applied for several applications. Formaggia, Gerbeau, Nobile and  Quarteroni \cite{formaggia1} have coupled 1-D and 3-D Navier-Stokes equations for studying blood flows in compliant vessels. In the context of  river dynamics, Miglio,  Perotto and  Saleri
 \cite{miglio}, Marin and Monnier  \cite{monnier}, Finaud-Guyot, Delenne, Guinot and  Llovel  \cite{finaud-guyot}, Malleron,  Zaoui, Goutal and  Morel
   \cite{malleron} have coupled 1-D and 2-D shallow water models. Leiva, Blanco and Buscaglia \cite{leiva2011} present also several such applications for Navier-Stokes equations. \par
Several techniques can be used to couple different models, either based on variational, algebraic, or domain decomposition approaches. In the context of dimensionally heterogeneous models, in addition to the previously mentioned references, there exists also a number of papers on this subject for purely hyperbolic problems (e.g. \cite{Godlewski2004,Godlewski2005}, \cite{bouttin}), but we will not elaborate on these studies since our focus is much more on hyperbolic/parabolic problems.
\par
In our study, we will focus on the design of an efficient Schwarz-like iterative coupling method.  The possibility of performing iterations between both models, i.e. of using a Schwarz method, is already considered in Miglio,  Perotto and  Saleri \cite{miglio} and Malleron, Zaoui, Goutal and Morel  \cite{malleron}. This kind of algorithm has several practical advantages. In particular, it is simple to develop and operate,  and it does not  require heavy changes in the numerical codes to be coupled: each model can be run separately,  the interaction between subdomains being ensured through boundary conditions only. These are important aspects in view of complex realistic applications.  
\par
Our final objective is to design an efficient  algorithm for the coupling of a 1-D/2-D shallow water model with a 2-D/3-D Navier-Stokes model. As a first step in this direction, the present study aims at identifying the main questions that we will have to face, as well as an adequate mathematical framework and possible ways to address these questions. We will perform this preliminary stage on a very simple testcase, coupling a 2-D Laplacian equation with a corresponding simplified 1-D equation. Seemingly similar testcases  were addressed by Blanco,
Discacciati and  Quarteroni \cite{blanco} and Leiva, Blanco and Buscaglia  \cite{leiva}, but with different  coupling methodologies
 (variational approach in \cite{blanco} and Dirichlet-Neumann coupling in \cite{leiva}). Moreover we have chosen to use non symmetrical boundary conditions in our 2-D model, in order to develop a fully two dimensional solution, and our 1-D model is obtained by integration of the 2-D equation along one direction, by analogy with the link between the shallow water system and the Navier-Stokes system. The rigorous mathematical derivation of the 1-D model clearly highlights its validity conditions.\par
Section 2 is devoted to the presentation of the 2-D Laplacian model, and to the derivation of the corresponding reduced 1-D model. Then a Schwarz iterative coupling algorithm is presented in Section 3, and its theoretical convergence properties are analyzed. In particular, the influence of the interface location is discussed. Finally numerical tests are presented in Section 4, which fully validate the previous analytical results.

\section{Derivation of the reduced model}\label{sec:derivation}

We are interested in the following boundary-valued problem in a domain $\w\subset\R^2$:
\begin{subequations}
\label{eq:full2-D0}
\begin{empheq}[left=\empheqlbrace\;]{align} 
&- \Delta u(x,z) = F(x,z),\quad\forall (x,z)\in \w,\label{eq:full2-D0a}\\
& \alpha\dn{u}(x,z)+\kappa u(x,z)= 0,\quad\forall (x,z)\in \partial\w,
\end{empheq}
\end{subequations}

with $\alpha$ and $\kappa$ are nonnegative numbers that allows for Dirichlet, Neumann or Robin boundary conditions.\\

\noindent In this section, we want to take advantage of the shallowness of the domain $\w$ (or some subdomain of $\w$) to derive a reduced model. As it is done for the derivation of the shallow water model (see \textit{e.g.} \cite{gerbeauperthame}), we want to replace the complete 2-D model \eqref{eq:full2-D0} by a (simpler) 1-D equation wherever it is possible (and keep the original 2-D model everywhere else). We will finally obtain the coupled 1-D / 2-D system \eqref{eq:1-Dmodel}-\eqref{eq:2-Dmodel} which we will analyze and simulate in the coming sections.\\
In order to discriminate between 1-D and 2-D regions, we introduce the following definition:
\begin{definition}\label{def1}
Let $\woned$ be the subset of $\w$ in which 2-D effects may be neglected, and $\wtwod=\w\backslash\woned$ the subset of $\w$ in which 2-D effects cannot be neglected.
\end{definition}
\begin{remark}\label{rem1}
Naturally the definition of $\woned$ depends on several features, such as the domain aspect ratio, the considered system of equations, forcing terms (including boudary conditions), etc.
\end{remark}

\begin{figure}[!ht]
\begin{center}
%\resizebox{12cm}{!}{\input{Figs/racket.pdf_t}}
\includegraphics[width=14cm]{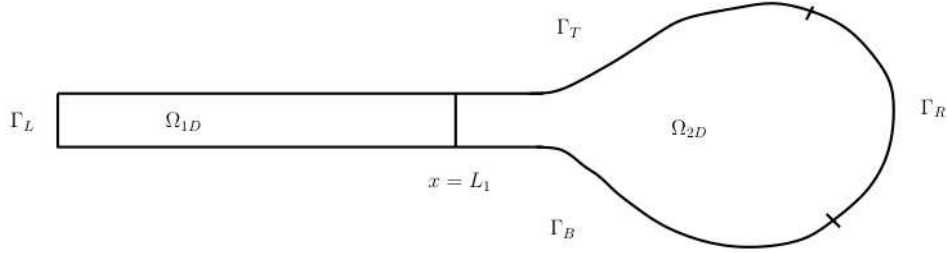}
\caption{\label{fig:racket}Typical computational domain $\w$, including a zone with a true 2-D behavior ($x\geq L_{1}$) together with a shallow zone where we intend to use a 1-D model ($x<L_{1}$).}
\end{center}
\end{figure}

\noindent For the sake of clarity (see Figure \ref{fig:racket}), we will assume that there exists $H$ and $L_1$ such that 
$$\woned=\w\cap\{x<L_{1}\}=(0,L_1) \times (0,H)\mbox{ and } \wtwod=\w\cap\{x> L_{1}\}. $$

%\begin{figure}[!ht]
%\begin{center}
%	\resizebox{12cm}{!}{\input{Figs/racket.pdf_t}}
%	\caption{\label{fig:racket}Typical computational domain $\w$, including a zone with a true 2-D behavior ($x\geq L_{1}$) together with a shallow zone where we intend to use a 1-D model ($x<L_{1}$).}
%\end{center}
%\end{figure}

Let us now consider equation \eqref{eq:full2-D0} with the following boundary conditions (see Figure \ref{fig:racket} for the notations):
\begin{subequations}
\label{eq:full2-D}
\begin{empheq}[left=\empheqlbrace\;]{align} 
- \Delta u(x,z) &= F(x,z),\quad\forall (x,z)\in \w,\\
\dn{u}(x,z)&= 0,\quad\forall (x,z)\in \Gamma_{T},\\
\dn{u}(x,z)+\kappa u(x,z)&= 0,\quad\forall (x,z)\in \Gamma_{B},\\
u(x,z) &= \gamma_1(x,z),\quad\forall (x,z)\in \Gamma_{L},\\
u(x,z) &= \gamma_2(x,z),\quad\forall (x,z)\in \Gamma_{R}.
\end{empheq}
\end{subequations}

In order to derive the 1-D model in $\woned=\w\cap\{x<L_1\}$, we introduce the following dimensionless variables and numbers:
\begin{eqnarray}
&\varepsilon=\frad{H}{L_{1}},\\
&\tilde{x} = \frad{x}{L_{1}},\quad\tilde{z} = \frad{z}{H},\quad \tilde{u}(\tilde{x}, \tilde{z})  
= \frad{u(x,z)}{U},\quad\tilde{F}(\tilde{x}, \tilde{z}) = F(x,z)\frad{ L_{1}^2}{U} \mbox{ and } \tilde{\kappa} = \kappa L_{1},
\end{eqnarray}
where $L_{1}$ (resp $H$) is the characteristic length (resp height) of $\woned$, $\varepsilon$ is called the aspect ratio, and $U$ is a characteristic value for $u(x,z)$.\\
The nondimensional form of equations \eqref{eq:full2-D} in $\woned$ reads\footnote{Since $\woned=(0,L_1)\times(0,H)$ we have $\vec n=\pm e_z$ in \eqref{eq:full2-Dndimb} and \eqref{eq:full2-Dndimc}.}:

\begin{subequations}
 \label{eq:full2-Dndim}
\begin{empheq}[left=\empheqlbrace\;]{align} 
&- \frad{\partial^2 \tilde{u}}{\partial \tilde{x}^2} - 
\frad{1}{\varepsilon^2} \frad{\partial^2 \tilde{u}}{\partial \tilde{z}^2} = \tilde{F}\quad 
\mbox{in}\quad \woned \label{eq:full2-Dndima} \\
&\frad{\partial \tilde{u}}{\partial \tilde{z}}= 0\quad \mbox{on}\quad \Gamma_{T}^{1} = \Gamma_T \cap \partial \woned \label{eq:full2-Dndimb}\\
&- \frad{1}{\varepsilon}\frad{\partial \tilde{u}}{\partial \tilde{z}} + \tilde{\kappa}  \tilde{u}= 0\quad \mbox{on}\quad \Gamma_{B}^{1}
= \Gamma_B \cap \partial \woned\label{eq:full2-Dndimc}\\[0.3cm]
&\tilde{u} = \frad{\gamma_1}{U}\quad \mbox{on}\quad \Gamma_L^{1} = \Gamma_L. %\;\;  u = \end{empheq}
\end{empheq}
\end{subequations}

%o\`u $\Gamma_h^1 = \Gamma_h \cap \partial \tilde{\Omega}_1$ et $\Gamma_f^2 = \Gamma_f \cap \partial \tilde{\Omega}_1$.\\
We assume (see \cite{REF-KAPPA} for the scaling of $\tilde\kappa$) that 
\be\label{eq:SW}
\tilde{F}= O(1), \frad{\partial^2 \tilde{u}}{\partial \tilde{x}^2} =O(1)\; \mbox{ and }\; \tilde{\kappa} = O(\varepsilon),
\ee
which is a sufficient condition to ensure that the 2-D effects are negligible in $\woned$. Indeed we deduce from equation \eqref{eq:full2-Dndima} that
\be
\frad{\partial^2 \tilde{u}}{\partial \tilde{z}^2} = O(\varepsilon^2).
\ee

By vertical integration on $(\tilde{z}, 1)$, and accounting for the boundary condition \eqref{eq:full2-Dndimb}, we find:
\be \label{eq:dzordre2}
\frad{\partial \tilde{u}}{\partial \tilde{z}} = O(\varepsilon^2)
\ee
and finally:
\be\label{eq:order2}
\tilde{u}(\tilde{x},\tilde{z}) = \tilde{u}(\tilde{x},0) + O(\varepsilon^2).
\ee
Going back to original variables, we have:
\begin{equation} \label{eqn : relation_ordre_2}
 u(x,z) = u(x,0) + O(\varepsilon^2), \;\forall \;(x,z) \in\woned.
\end{equation}

\begin{figure}[h]
\begin{center}
\includegraphics[width=8cm]{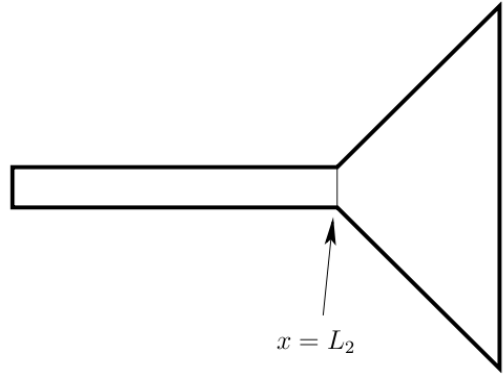}
	\caption{\label{fig:entonnoir}Funnel-shaped computational domain. The domain is shallow for $x<L_{2}$.}
\end{center}
\end{figure}

\begin{figure}[h]
\begin{center}
	\includegraphics[width=15cm]{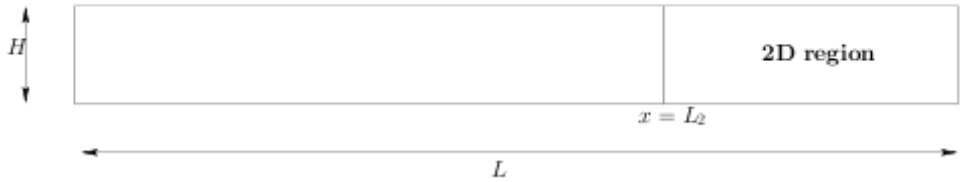}	
	\caption{\label{fig:rectangle}Rectangular computational domain $\w$. The domain is shallow: $H/L\ll 1$, and we assume that the forcing terms are supported in $\{x>L_{2}\}$.}
\end{center}
\end{figure}

We now introduce the averaging operator in the vertical direction. For any function $f$ of $z$, we set:
\be
\av{f}=\frad{1}{H}\int_{0}^Hf(z)\,dz.
\ee
We integrate equation \eqref{eq:order2} for $z \in (0, H)$ and obtain:
\begin{equation} \label{eqn : relation_ordre_2_bis}
 \bar{u}(x) = u(x,0) + O(\varepsilon^2),\;\forall x \in [0, L].
\end{equation}
We now average equation \eqref{eq:full2-Dndima} in the vertical direction, taking into account the Robin boundary condition \eqref{eq:full2-Dndimc} on $\Gamma_B^{1}$, and find:
\be
-\frad{\partial^2 \bar{u}}{\partial x^2}  + \frad{\kappa}{H} u(x,0) = \bar{F}, 
\ee
%o\`u $\displaystyle \bar{F} = \frad{1}{\varepsilon} \int_{0}^{\varepsilon} F(x,z) dz$.\\
For every $x \in [0, L_{1}]$ we may use approximation (\ref{eqn : relation_ordre_2_bis}) to introduce the new problem:
\begin{equation} \label{eqn : chaleur1-D_robin}
-\frad{\partial^2  u_1}{\partial x^2} + \frad{\kappa}{H} u_1 = \bar{F} \quad \mbox{in}\quad [0, L_{1}].
\end{equation}
It will replace \eqref{eq:full2-D0a} in $\woned$. As evoked in Remark \ref{rem1}, the reader can be easily convinced that it is particularly awkward to guess the value of $L_{1}$. Indeed one has to specify the criteria that define 2-D effects, and in practical situations we may only be able to define $L_{2}$ which is such that $\big(\w\cap\{x\geq L_{2}\}\big)\subset\wtwod$, or in other words $L_{2}\geq L_{1}$. In this work we consider two different situations:
\begin{list}{-}{}
\item a funnel-shaped domain (see Figure \ref{fig:entonnoir}) with a thin left part, so that we anticipate 2-D effects on the right (wide) part of the domain. In this case the definition of $L_{2}$ is based on a geometrical criterion.
\item a rectangular domain (see Figure \ref{fig:rectangle}) with a small aspect ratio $\varepsilon=H/L$ (so that we can anticipate weak 2-D effects), but with some 2-D forcing terms occuring in the right end of the domain. In that case the definition of $L_{2}$ is based on the support of forcing terms.
\end{list}

%\begin{figure}[h]
%\begin{center}
%	\resizebox{6cm}{!}{\input Figs/Entonnoir.pdf_t}
%	\caption{\label{fig:entonnoir}Funnel-shaped computational domain. The domain is shallow for $x<L_{2}$.}
%\end{center}
%\end{figure}
%
%\begin{figure}[h]
%\begin{center}
%	\resizebox{15cm}{!}{\input Figs/Rectangle.pdf_t}
%	\caption{\label{fig:rectangle}Rectangular computational domain $\w$. The domain is shallow: $H/L\ll 1$, and we assume that the forcing terms are supported in $\{x>L_{2}\}$.}
%\end{center}
%\end{figure}

\noindent At this point we have defined an upper bound $L_{2}\geq L_{1}$, but the exact value of $L_{1}$ remains unknown. From now on we choose an interface $L_{0}$ without any \textit{a priori} information (other than $L_{0}<L_{2}$) and decide to consider the model reduction \eqref{eq:1-Dmodel} on $\w_{1}=\w\cap\{x<L_{0}\}$, while we keep the 2-D model in $\w_{2}=\w\cap\{x>L_{0}\}$. Finally we have the two following systems
\begin{eqnarray} \label{eq:1-Dmodel}
\mbox{1-D model: } \qquad\left\{
\begin{array}{rcl}
\displaystyle -\frad{\partial^2 u_1}{\partial x^2} + \frad{\kappa}{H} u_1 &=& F_1  \mbox{ in } (0, L_0),\\
u_1(0)&=& \bar{\gamma}_1.
\end{array}
\right.
\end{eqnarray}
and
\begin{eqnarray} \label{eq:2-Dmodel}
\mbox{2-D model: } \qquad \left\{
\begin{array}{rcl}
\displaystyle  - \Delta u_2 &=& F_2 \mbox{ in } \w_2,\\
\displaystyle \frad{\partial u_2}{\partial n}= 0 \mbox{ on } \Gamma^2_{T} &=& \Gamma_T \cap \partial \w_2, \\
\displaystyle  \frad{\partial u_2}{\partial n} + \kappa u_2= 0 \mbox{ on }\Gamma^2_{B} &=& \Gamma_B \cap \partial \w_2,  \\
 u_2 &=& \gamma_2 \mbox{ on } \Gamma_R.
\end{array}
\right.
\end{eqnarray}
where $F_1 = \overline{F}$ and $F_2 = F_{|\Omega_2}$.\\

Two cases may occur:
\begin{list}{-}{}
\item Favourable case: $L_{0}<L_{1}$, so that $\w_{1}\subset\woned$ and the following model reduction is relevant. In particular, hypothesis \eqref{eq:SW} holds so that Theorem \ref{th:errorcontrol} applies,
\item Unfavourable case: $L_{0}\geq L_{1}$, so that $\w_{1}\not\subset\woned$ and the 1-D model will not be able to reproduce the 2-D reality (in particular, hypothesis \eqref{eq:SW} does not hold).
\end{list}

We now want to evaluate this model reduction in the favourable case $L_{0}<L_{1}$.

\section{Coupling algorithm}\label{algorithm}
\begin{figure}
\begin{center}
\includegraphics[width=14cm]{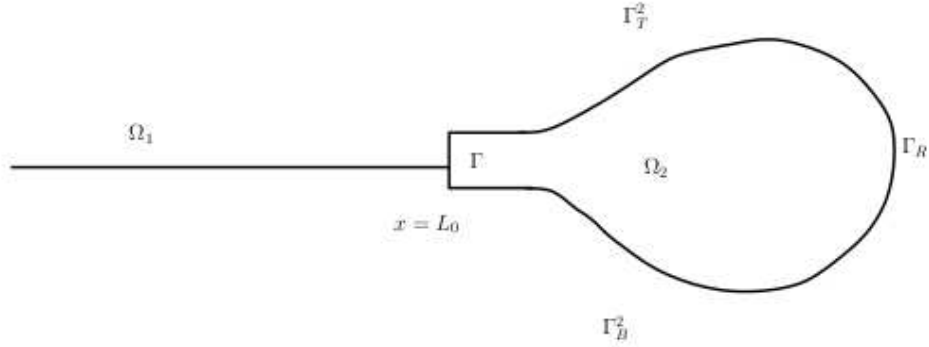}
 \caption{Computational domains for the 1-D/2-D reduced model}
\label{fig : domaine_couplage_1-D_2-D}
 \end{center}
\end{figure}

Let us consider the two models (\ref{eq:1-Dmodel}) and (\ref{eq:2-Dmodel}) to be coupled respectively through the interfaces $x = L_0$ and $\Gamma$ as shown in Figure \ref{fig : domaine_couplage_1-D_2-D}. \\
%\begin{figure}
%\begin{center}
%  \scalebox{0.7}{\input{Figs/racket_couplage.pdf_t}}
% \caption{Computational domains for the 1-D/2-D reduced model}
%\label{fig : domaine_couplage_1-D_2-D}
% \end{center}
%\end{figure}
In coupling problems, the first difficulty lies in defining the coupling notion by itself, \textit{i.e.} defining the quantities or values to be exchanged between the two 
models through the coupling interfaces. In our case and from a physical point of view, one may propose the following conditions, see \cite{blanco}, \cite{leiva}:
\begin{subequations}
\label{contraintes}
\begin{empheq}[left=\empheqlbrace\;]{align} 
u_1(L_0) &= \frac{1}{H} \int_{0}^{H} u_2(L_0, z) dz  \label{eqn : contrainte1} \\
\frac{\partial u_1}{\partial x}(L_0) &=\frac{1}{H} \int_{0}^{H} \frac{\partial u_2}{\partial x}(L_0, z) dz \label{eqn : contrainte2}
\end{empheq}
\end{subequations}
which correspond to the conservation of $u$ and its flux through the interface.\\
Unfortunately these two constraints do not allow the well-posedness of the 2-D model and of the coupled problem. They are called \textit{defective boundary conditions} in the 
literature, \cite{formaggia1}, \cite{formaggia2}, \cite{leiva}.\\
One should then rather apply a coupling method ensuring the following points:
\begin{itemize}
\item [(i)] the well-posedness of the 1-D and 2-D models
\item [(ii)] the physical constraints are satisfied
\item [(iii)] the control of the difference between the coupled solution and the reference one (corresponding to the 2-D model over the whole domain $\w$). Indeed,
due to the nature of the problem, the reader can be easily convinced that one does not expect to end with  a solution of the 2-D model
equal to the restriction of the reference one on $\Omega_2$. 
\end{itemize}
The coupling problem with (\ref{eqn : contrainte1}) and (\ref{eqn : contrainte2}) has been studied in \cite{blanco} and \cite{leiva} using variational and algebraic approaches. 
%Both in these two papers,
%authors propose a constant extension of (\ref{eqn : contrainte1}) or (\ref{eqn : contrainte2}).\\
In this section, we propose an iterative coupling method based on classical Schwarz algorithms. These iterative methods were used for the first time in
the context of dimensionally heterogeneous coupling in \cite{formaggia1} and \cite{miglio} to study
a nonlinear hyperbolic coupling problem.\\
We will prove the convergence of these algorithms given an appropriate choice of boundary conditions  at $x=L_0$ and on $\Gamma$. 
Then we will study the solutions obtained after convergence
and compare them to the global reference  solution $u$ defined by (\ref{eq:full2-D}). Finally we will give some results regarding the choice of the coupling interface position.\\ 

\subsection{Schwarz algorithms}
Let us introduce first the Schwarz algorithms in the context of dimensionally homogeneous coupling.
Consider the two systems, defined on $\Omega_1$ and $\Omega_2$ shown in Figure \ref{fig : domaine_couplage}:
\begin{figure}
\begin{center}
\includegraphics{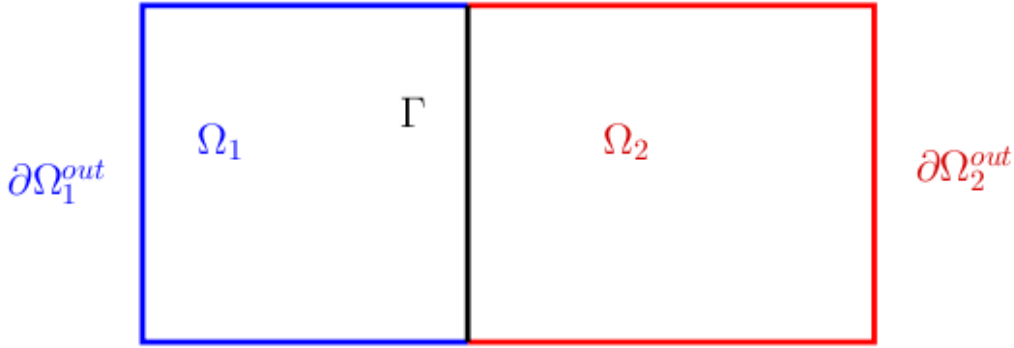}
 \caption{Computational domains for the dimensionally homogeneous  coupling problem}
  \label{fig : domaine_couplage}
 \end{center}
\end{figure}
%\begin{figure}
%\begin{center}
% \input{Figs/DDM.pdf_t} 
% \caption{Computational domains for the dimensionally homogeneous  coupling problem}
%  \label{fig : domaine_couplage}
% \end{center}
%\end{figure}
\begin{eqnarray}\label{eqn : modele_mono}
\left\{
\begin{array}{ll}
 \mathcal{L}_u u = f_u \;\;\; \mbox{in}\;\;\; \Omega_1 \subset \mathbb{R}^n\\
B_u^{out} u = g_u \;\;\; \mbox{on} \;\;\; \partial \Omega_1^{out} \\
\end{array}
\right.
\;\;\;\mbox{and}\;\;\;\;\;
\left\{
\begin{array}{ll}
 \mathcal{L}_v v = f_v \;\;\; \mbox{in}\;\;\; \Omega_2 \subset \mathbb{R}^n \\
B_v^{out} v = g_v \;\;\; \mbox{on} \;\;\; \partial \Omega_2^{out} \\
\end{array}
\right.
\end{eqnarray}
The operators $\mathcal{L}_u$ and $\mathcal{L}_v$ are different. We assume that $u$ and $v$ have to satisfy  the following constraints derived from the physics :
\begin{subequations}
\label{eq:mono}
\begin{empheq}[left=\empheqlbrace\;]{align} 
C_1 u &= C_2 v \label{eqn : condition_couplage1}\\
C'_1 u &= C'_2 v \label{eqn : condition_couplage2}
\end{empheq}
\end{subequations}
through the interface $\Gamma$, where $C_1$, $C'_1$, $C_2$ and $C'_2$ are differential operators.\\
To couple these two models we can implement the following iterative algorithm:\\
{\it For a given $ v^0 $  and at each iteration $k \geq 0$, solve:
$$
\left\{
\begin{array}{lll}
 \mathcal{L}_u u^{k+1} &=& f_u \quad \mbox{in}\quad \Omega_1 \\
B_u^{out} u^{k+1} &=& g_u \quad \mbox{on} \quad \partial \Omega_1^{out} \\
B_u u^{k+1} &=& B_v v^k \quad \mbox{on}\quad \Gamma
\end{array}
\right.
\quad \mbox{then}\quad
\left\{
\begin{array}{lll}
 \mathcal{L}_v v^{k+1} &=& f_v \quad \mbox{in}\quad \Omega_2 \\
B_v^{out} v^{k+1} &=& g_v \quad \mbox{on} \quad \partial \Omega_2^{out} \\
B'_v v^{k+1} & =& B'_u u^{k+1} \quad \mbox{on}\quad \Gamma.
\end{array}
\right.
$$}
Once convergence is achieved, the physical constraints (\ref{eqn : condition_couplage1}) and (\ref{eqn : condition_couplage2}) have to be satisfied. Then
care should be taken to choose  
the operators $B_u$, $B'_u$, $B_v$ and $B'_v$ in order to ensure convergence toward the unique solution
defined by (\ref{eqn : modele_mono}),  (\ref{eqn : condition_couplage1}) and 
(\ref{eqn : condition_couplage2}), see \cite{quarteroni} and \cite{vero}. \\
This method can be generalized to the dimensionally heterogeneous coupling case.
Let assume that we have to solve the following 1-D model/2-D model coupled problem:
$$
\left\{
\begin{array}{ll}
 \mathcal{L}_1 u_1 = f_1 \quad \mbox{in}\quad \Omega_1 \subset \mathbb{R}\\
B_1^{out} u_1 = g_1 \quad \mbox{on} \quad \partial \Omega_1^{out} \\
\end{array}
\right.
\quad\mbox{and}\quad
\left\{
\begin{array}{ll}
 \mathcal{L}_2 u_2 = f_2 \quad \mbox{in}\quad \Omega_2 \subset \mathbb{R}^2\\
B_2^{out} u_2 = g_2 \quad \mbox{on} \quad \partial \Omega_2^{out} \\
\end{array}
\right.
$$
and  we suppose that we have the following coupling constraints to satisfy at $x=L_0$ and on $\Gamma$:
\begin{subequations}
\label{mono}
\begin{empheq}[left=\empheqlbrace\;]{align} 
C_1 u_1(L_0) &= C_2 \left(\mathcal{R}u_2\right) (L_0)  \label{eq : contrainte_generale1}\\
C'_2 u_2(L_0,z) &= C'_1 \left(\mathcal{E}u_1\right) (L_0,z)\quad\mbox{on}\quad \Gamma \label{eq : contrainte_generale2}
\end{empheq}
\end{subequations}
$\mathcal{R}u_2$ is a  restriction of $u_2$ at $x=L_0$ and $\mathcal{E}u_1$ is an extension of $u_1(L_0)$ all along $\Gamma$. More generally
we can define the operators $\mathcal{R}$ and $\mathcal{E}$ as in \cite{blanco} by:
\begin{eqnarray*}
 \mathcal{R} : \Lambda& \longrightarrow &\Lambda_0 \\
u_{2|\Gamma} &\longmapsto& \mathcal{R}u_{2|\Gamma}
\end{eqnarray*}
and
\begin{eqnarray*}
 \mathcal{E} : \Lambda_0& \longrightarrow &\Lambda \\
u_{1|x=L_0} &\longmapsto& \mathcal{E}u_{1|x = L_0}
\end{eqnarray*}
The spaces $\Lambda_0$ et $\Lambda$ are the trace spaces on the interface $x=L_0$ for 1-D functions and on the interface $\Gamma$
for 2-D functions. As mentioned in \cite{blanco}, these two operators are not invertible.\\
 One may thus
implement the following algorithm:\\
{\it For a given $ u_2^0 $  and at each iteration $k \geq 0$, solve}:
$$
\left\{
\begin{array}{lll}
 \mathcal{L}_1 u_1^{k+1} &=& f_1 \quad \mbox{in}\quad \Omega_1 \\
B_1^{out} u_1^{k+1} &=& g_1 \quad \mbox{on} \quad \partial \Omega_1^{out} \\
B_1 u_1^{k+1} &=& B_2 \mathcal{R}u_2^k \quad \mbox{at}\quad x=L_0
\end{array}
\right.
\quad then \quad
\left\{
\begin{array}{lll}
 \mathcal{L}_2 u_2^{k+1} &=& f_2 \quad \mbox{in}\quad \Omega_2 \\
B_2^{out} u_2^{k+1} &=& g_2 \quad \mbox{on} \quad \partial \Omega_2^{out} \\
B'_2 u_2^{k+1} & =& B'_1 \mathcal{E}u_1^{k+1} \quad \mbox{on}\quad \Gamma.
\end{array}
\right.
$$
In practice, we do not have conditions such as (\ref{eq : contrainte_generale1}) (at $x=L_0$) and (\ref{eq : contrainte_generale2}) (for $(x,z)\mbox{ on }\Gamma$), but only conditions at $x=L_0$ such as (\ref{eqn : contrainte1}) and (\ref{eqn : contrainte2}). 
This leads us to make a choice of the operators $\mathcal{R}$, $\mathcal{E}$, $C_1$, $C_2$, $C'_1$ and $C'_2$.\\
  In general the choice of the restriction operator $\mathcal{R}$ is more straightforward than the choice of the extension operator $\mathcal{E}$ one. In this
study, since the 1-D model is obtained after some approximations and by averaging the 2-D model, it is reasonable to define $\mathcal{R}$ as  the vertical average. On the other hand, the question of the choice of the operator $\mathcal{E}$ remains open. \\
In \cite{blanco} and \cite{leiva}, authors proposed a constant extension of (\ref{eqn : contrainte1}) and (\ref{eqn : contrainte2}) along $\Gamma$, and 
impose the following strong coupling constraints:
$$
\left\{
\begin{array}{lll}
 u_2(L_0, z) &= & u_1(L_0)\quad \mbox{on} \quad \Gamma\\
\displaystyle \frac{\partial u_1}{\partial x}(L_0) &=& \displaystyle \frac{1}{H} \int_{0}^{H} \frac{\partial u_2}{\partial x}(L_0, z) dz
\end{array}
\right.
\;\;\;\mbox{or}\;\;\;\;\;
\left\{
\begin{array}{lll}
 u_1(L_0) &= & \displaystyle \frac{1}{H}\int_{0}^{H}u_2(L_0,z) dz\\
\displaystyle \frac{\partial u_2}{\partial x}(L_0,z) &=& \displaystyle \frac{\partial u_1}{\partial x}(L_0)\quad \mbox{on}\quad \Gamma
\end{array}
\right.
$$
It is a choice among many others. One  may also choose a multitude of operators $C_1$, $C_2$, $C'_1$ and $C'_2$ ensuring the physical constraints to be satisfied.\\
The strategy that we adopt here is to
choose, due to  relations (\ref{eqn : relation_ordre_2}) and (\ref{eqn : relation_ordre_2_bis}), a constant extension of $u_1$ along $\Gamma$ 
and then to implement a family of Schwarz algorithms with appropriate boundary 
conditions at $x = L_0$ and on $\Gamma$. 
In this case Schwarz coupling algorithm reads:\\
{\it For a given $ u_2^0 $  and at each iteration $k \geq 0$, solve :
\begin{eqnarray} \label{eqn : iteration_1-D}
\left\{
\begin{array}{lllll}
\displaystyle - \dxx{u_1^{k+1}} + \frac{\kappa}{H} u_1^{k+1}= F_1 \quad \mbox{in}\quad (0, L_0) \\
\\
u_1^{k+1}(0) = \bar{\gamma}_1  \\
\\
\displaystyle B_1 u_1^{k+1}(L_0) = B_1\bar{u}_2^k
\end{array}
\right.
\end{eqnarray}
and then solve
\begin{eqnarray} \label{eqn : iteration_2-D}
\left\{
\begin{array}{lllllll}
-\Delta u_2^{k+1} = F_2 \quad \mbox{in}\quad \wtwo \\
\\
 \displaystyle \dn{u_2^{k+1}}= 0 \quad \mbox{on}\quad \Gamma_{T}^2 \\
\\
\displaystyle  \dn{u_2^{k+1}}+ \kappa u_2^{k+1}= 0 \quad \mbox{on}\quad \Gamma_{B}^2 \\
\\
u_2^{k+1} = \gamma_2 \quad \mbox{on} \quad \Gamma_R\\
\\
\displaystyle B_2 u_2^{k+1} = B_2 u_1^{k+1} \quad \mbox{on}\quad \Gamma
\end{array}
\right.
\end{eqnarray}
}
The linear operators $B_1$ and $B_2$ will be defined such that the points (i), (ii) and (iii) (see introduction of Section \ref{algorithm}) are satisfied and such that the algorithm converges.\\
We will first study the convergence of the coupling algorithm. Subsequently we move on the point (iii).\\
To ensure the convergence of Schwarz algorithms in the case of classical domain decomposition without overlapping, it is proposed in \cite{lions}
to use Robin  operators. We will extend the use of these operators to our coupling problem. We define the operators $B_1$ and $B_2$ for a given $\lambda > 0$ as follows:
\begin{eqnarray} \label{eqn : robin_operator1}
B_1= \frac{\partial}{\partial n_1} + \lambda \id
\end{eqnarray}
and
\begin{eqnarray} \label{eqn : robin_operator2}
B_2 = \dntwo{} + \lambda \id
\end{eqnarray}
where $n_1$ and $n_2$ are the outward unit normal to the 1-D and 2-D domains respectively.
% \textcolor{red}{These operators are also called {\it one-sided} Robin operators. }\\
We note that the operators $B_1$ and $B_2$ ensure the well-posedness of the problem at each iteration.  Let us study the convergence of Schwarz 
algorithm with this family of operators.
\subsection{Algorithm convergence}
\begin{prop}
 For each $ \lambda > 0$, the $(u_1^k, u_2^k)$ Schwarz algorithm converges in  $H^1(\Omega_1) \times H^1(\Omega_2)$ 
 to $(u_1^{\lambda}, u_2^{\lambda})$  
 that satisfies the physical constraints (\ref{eqn : contrainte1}), 
 (\ref{eqn : contrainte2}).\\
\end{prop}
\noindent {\bf Proof:}\\
Let define the differences between two successive iterations:
 $$
 e_1^{k+1}(x) = u_1^{k+1}(x) - u_1^k(x), \quad\forall x \in (0, L_0)
 $$
 and
 $$
 e_2^{k+1}(x,z) = u_2^{k+1}(x,z) - u_2^k(x,z),\quad \forall (x,z)  \in \wtwo
 $$
These functions satisfy the following systems:
\begin{eqnarray} \label{eqn : error_1-D}
\left\{
\begin{array}{lllll}
\displaystyle - \dxx{e_1^{k+1}} + \frac{\kappa}{H} e_1^{k+1}= 0 \quad \mbox{in}\quad (0, L_0) \\
\\
e_1^{k+1}(0) = 0  \\
\\
\displaystyle  \frac{\partial e_1^{k+1}}{\partial x}(L_0) + \lambda e_1^{k+1}(L_0) = 
\frac{\partial \bar{e}_2^{k}}{\partial x}(L_0) + \lambda \bar{e}_2^{k}(L_0)
\end{array}
\right.
\end{eqnarray}
and
\begin{eqnarray} \label{eqn : error_2-D}
\left\{
\begin{array}{lllllll}
-\Delta e_2^{k+1} = 0 \quad \mbox{in}\quad \wtwo \\
\\
 \displaystyle \dn{e_2^{k+1}}= 0 \quad \mbox{on}\quad \Gamma_{T}^2 \\
\\
\displaystyle  \dn{e_2^{k+1}}+ \kappa e_2^{k+1}= 0 \quad \mbox{on}\quad \Gamma_{B}^2 \\
\\
e_2^{k+1} = 0 \quad \mbox{on} \quad \Gamma_R\\
\\
\displaystyle  -\frac{\partial e_2^{k+1}}{\partial x}(L_0,z)  + \lambda e_2^{k+1}(L_0,z)
= -\frac{\partial e_1^{k+1}}{\partial x}(L_0)  + \lambda e_1^{k+1}(L_0)\quad \mbox{on}\quad \Gamma.
\end{array}
\right.
\end{eqnarray}
The first two equations of (\ref{eqn : error_1-D}) lead to:
\begin{equation} \label{eqn : error1D}
e_1^{k+1}(x) = \alpha_{k+1} \sinh(a x), \;\;\; \forall x \in (0, L_0)
\end{equation}
where $\alpha_{k+1} \in \mathbb{R}$ and $\displaystyle a = \sqrt{\frac{\kappa}{H}}$.\\
If we take the vertical average of the boundary condition on $\Gamma$ in (\ref{eqn : error_2-D}), and due to the boundary condition at $x=L_0$, we obtain :
$$
\left\{
\begin{array}{ccc}
\displaystyle -\frac{\partial \bar{e}_2^{k+1}}{\partial x}(L_0)  + \lambda \bar{e}_2^{k+1}(L_0)&= &
\displaystyle-\frac{\partial e_1^{k+1}}{\partial x}(L_0)  + \lambda e_1^{k+1}(L_0)\\
\\
\displaystyle \frac{\partial \bar{e}_2^{k+1}}{\partial x}(L_0) + \lambda \bar{e}_2^{k+1}(L_0)& = &
\displaystyle \frac{\partial e_1^{k+2}}{\partial x}(L_0) + \lambda e_1^{k+2}(L_0).
\end{array}
\right.
$$
This implies:
\begin{eqnarray} \label{eqn : moyenne_e2}
 \bar{e}_2^{k+1}(L_0) &= & \frac{1}{2\lambda}\left(A \alpha_{k+2} + B \alpha_{k+1} \right)
\end{eqnarray}
and
\begin{eqnarray} \label{eqn : moyenne_derivee_e2}
\dx{\bar{e}_2^{k+1}}(L_0)  = \frac{1}{2} \left( A \alpha_{k+2} - B \alpha_{k+1}\right)
\end{eqnarray}
where $A = a \cosh(a L_0) + \lambda \sinh(a L_0)$ and $B = -a \cosh(a L_0) + \lambda \sinh(a L_0)$.\\
Now by multiplying the first equation of (\ref{eqn : error_2-D}) by $e_2^{k+1}$ and by integrating in $\Omega_2$, we obtain:
\begin{eqnarray*}
\int_{\Omega_2}^{} |\nabla e_2^{k+1}|^2 dx dz - \int_{\partial \Omega_2}^{} \frac{\partial e_2^{k+1}}{\partial n} e_2^{k+1}  d \sigma = 0
\end{eqnarray*}
then using the boundary conditions on $\Gamma_T^2$ and $\Gamma_R$, we obtain:
\begin{eqnarray} \label{eqn : FV_e2}
\int_{\Omega_2}^{} |\nabla e_2^{k+1}|^2 dx dz + \int_{\Gamma_B^2}^{} \kappa |e_2^{k+1} |^2 dx  &= &
-\int_{\Gamma}^{} \frac{\partial e_2^{k+1}}{\partial x}(L_0,z) e_2^{k+1}(L_0,z)  dz.
\end{eqnarray}
We replace $\displaystyle \frac{\partial e_2^{k+1}}{\partial x}(L_0,z)  $ in (\ref{eqn : FV_e2}) by its value obtained by using Robin boundary condition on $\Gamma$:
\begin{eqnarray} \label{eqn : FV_e2bis}
\int_{\Omega_2}^{} |\nabla e_2^{k+1}|^2 dxdz + \int_{\Gamma_B^2}^{} \kappa |e_2^{k+1}|^2 dx  &= & \int_{\Gamma}^{} 
\left( - \frac{\partial e_1^{k+1}}{\partial x}(L_0) + \lambda e_1^{k+1}(L_0) - \lambda e_2^{k+1}(L_0,z)\right) e_2^{k+1}(L_0,z) dz \nonumber\\
&=& B \alpha_{k+1} H \bar{e}_2^{k+1}(L_0) - \lambda  \int_{0}^{H}|e_2^{k+1} |^2(L_0,z) d z \nonumber\\
& = & \frac{B \alpha_{k+1} H}{2 \lambda}\left(A \alpha_{k+2} + B \alpha_{k+1} \right) -
 \lambda  \int_{0}^{H}|e_2^{k+1} |^2(L_0,z) d z.
\end{eqnarray}
We now  replace $e_2^{k+1}(L_0,z)$ in (\ref{eqn : FV_e2}) using the same Robin boundary condition:
\begin{eqnarray*}
\int_{\Omega_2}^{} |\nabla e_2^{k+1}|^2 dxdz + \int_{\Gamma_B^2}^{} \kappa |e_2^{k+1} |^2 d \sigma  &= & -\int_{\Gamma}^{} 
\frac{1}{\lambda}\left( - \frac{\partial e_1^{k+1}}{\partial x}(L_0) + \lambda e_1^{k+1}(L_0)\right)\frac{\partial e_2^{k+1}}{\partial x}(L_0,z) dz\\
&& - \frac{1}{\lambda}
 \int_{0}^{H}\left|\frac{\partial e_2^{k+1}}{\partial x}(L_0,z) \right|^2 d z \\
&=& - \frac{B \alpha_{k+1}H}{\lambda} \frac{\partial \bar{e}_2^{k+1}}{\partial x}(L_0) - \frac{1}{\lambda}\int_{0}^{H}\left|
\frac{\partial e_2^{k+1}}{\partial x}(L_0,z) \right|^2 d z \\
& = & -\frac{B \alpha_{k+1} H}{2 \lambda}\left(A \alpha_{k+2} - B \alpha_{k+1} \right) - \frac{1}{\lambda}
 \int_{0}^{H}\left|\frac{\partial e_2^{k+1}}{\partial x}(L_0,z) \right|^2 d z.
\end{eqnarray*}
Due to the fact that $\lambda > 0$, we deduce that:
$$
\frac{B \alpha_{k+1} H}{2 \lambda}\left(A \alpha_{k+2} + B \alpha_{k+1} \right) \geq 0
$$
and
$$
-\frac{B \alpha_{k+1} H}{2 \lambda}\left(A \alpha_{k+2} - B \alpha_{k+1} \right) \geq 0.
$$
Thus:
$$
A^2 \alpha_{k+2}^2 - B^2 \alpha_{k+1}^2 \leq 0.
$$
Then we obtain: 
$$
\frac{\alpha_{k+2}^2}{\alpha_{k+1}^2} \leq \frac{B^2}{A^2} = 
\left|\frac{-a \cosh(a L_0) + \lambda \sinh(a L_0)}{a \cosh(a L_0) + \lambda \sinh(a L_0)} \right|^2 < 1
$$
and finally
\begin{equation} \label{eqn : suite_geo}
\left| \frac{\alpha_{k+2}}{\alpha_{k+1}} \right| <  \left| \frac{B}{A} \right| < 1.
\end{equation}
So that the sequence $(\alpha_k)_{k \in \mathbb{N}}$ converge to zero.\\
%Now we will prove that $(u_1^k)_{k \in \mathbb{N}}$ and $(u_2^k)_{k \in \mathbb{N}}$ are Cauchy sequences respectively in $H^1(\Omega_1)$  and $H^1(\Omega_2)$.\\
Let us now remark that for all $k \geq 0$, $n \geq 0$, we have:
\begin{eqnarray*}
 u_1^{k+n} - u_1^{k} &=& \sum_{p=0}^{n-1} e_1^{k+p+1}\\
\end{eqnarray*}
and
\begin{eqnarray*}
 \frac{ \partial \left( u_1^{k+n} - u_1^{k} \right)}{\partial x} &=&  \sum_{p=0}^{n-1} \frac{\partial  e_1^{k+p+1} }{\partial x}\\
\end{eqnarray*}
Using the relation (\ref{eqn : error1D}) and the fact that the sequence $(\alpha_k)_{k \in \mathbb{N}}$ converges, we can
prove that $(u_1^k)_{k \in \mathbb{N}}$ and $\displaystyle (\frac{\partial u_1^k}{\partial x})_{k \in \mathbb{N}}$ are Cauchy sequences in
$L^2(\Omega_1)$. So that $(u_1^k)_{k \in \mathbb{N}}$ is a Cauchy sequence in $H^1(\Omega_1)$.\\
In the same way we observe that for all $k \geq 0$, $n \geq 0$, we have:
\begin{eqnarray*}
\nabla (u_2^{k+n} - u_2^{k}) &=& \sum_{p=0}^{n-1} \nabla e_2^{k+p+1}\\
\end{eqnarray*}
Using (\ref{eqn : FV_e2bis}), we deduce that:
$$
\int_{\Omega_2}^{} |\nabla e_2^{k+1}|^2 dxdz \leq  \frac{B \alpha_{k+1} H}{2 \lambda}\left(A \alpha_{k+2} + B \alpha_{k+1} \right)
$$
and then we can prove that $(\nabla e_2^k)_{k \in \mathbb{N}}$ is a Cauchy sequence in $L^2(\Omega_2)$ and due to the Poincar\'e inequality we have also
 $(u_2^k)_{k \in \mathbb{N}}$ is a Cauchy sequence in $L^2(\Omega_2)$. So that $(u_2^k)_{k \in \mathbb{N}}$ is a Cauchy sequence in $H^1(\Omega_2)$.\\
To conclude we have prove that $(u_1^k, u_2^k)$ Schwarz algorithm converges in  $H^1(\Omega_1) \times H^1(\Omega_2)$. 
 Moreover, at convergence the limit $(u_1^{\lambda}, u_2^{\lambda})$ verifies $B_2 u_2^{\lambda} = B_2 u_1^{\lambda}$
      and $B_1 u_1^{\lambda} = B_1 \bar{u}_2^{\lambda}$. Taking the 
vertical average on $\Gamma$ gives two linear combinations of the constraints  (\ref{eqn : contrainte1}), 
 (\ref{eqn : contrainte2}).$\Box$

\begin{remark}~\\
\begin{itemize}
\item  The Schwarz algorithms converge for all  $\lambda$ positive, but we remark that for $\displaystyle \lambda = a \coth(a  L_0)$, we have exact
convergence in two iterations. Indeed we have in this case:
$$
\displaystyle -\frac{\partial e_1^{k+1}}{\partial x}  + a \coth(a L_0) e_1^{k+1} =0 \;\;\; \forall k \geq 0,
$$
and then:
$$
B_2 e_2^{k+1} = -\frac{\partial e_2^{k+1}}{\partial x}  + a \coth(a L_0) e_2^{k+1} =-\frac{\partial e_1^{k+1}}{\partial x}  + a \coth(a L_0) e_1^{k+1}=0 \;\;\; \forall k \geq 0,
$$
The operator $\displaystyle \dnone{} + a \coth(aL_0) \id$ corresponds to the {\it absorbing} operator of the 1-D model.\\
We denote by $\lambda_{opt}$ this value of $\lambda$, see for example \cite{japhetnataf}.
\item If we take $\lambda_1 \neq \lambda_2$, we have {\it a priori} $(u_1^{\lambda_1}, u_2^{\lambda_1}) \neq (u_1^{\lambda_2}, u_2^{\lambda_2})$. 
This is in accordance with the ill-posedness of coupling problem defined by (\ref{eq:1-Dmodel}), (\ref{eq:2-Dmodel}), (\ref{eqn : contrainte1}) and (\ref{eqn : contrainte2}).
\end{itemize}
For the sake of clarity we will denote $(u_1, u_2)$ the limit of Schwarz algorithm instead of $(u_1^{\lambda}, u_2^{\lambda})$.
\end{remark}

\subsection{Control of the difference between the coupled solution and the global reference solution}
Unlike the case of domain decomposition, at convergence of Schwarz algorithm, we have $u_2 \neq u|_{\Omega_2}$ due to the model reduction. But as  mentioned above, we have chosen the family of {\it Robin} operators in order to get some control of
 the difference between $u_2$ and $u|_{\Omega_2}$. In fact we have the following result:
\begin{theorem}\label{th:errorcontrol}
  for each $\lambda > 0$, let  $(u_1^{\lambda}, u_2^{\lambda})$ denotes the limit of the Schwarz algorithm. If $L_0 < L_1$ then 
there exists $M(\lambda) > 0$ such that:
\begin{eqnarray} \label{eq : error_control}
\| u_{|\Omega_2} - u_2 \|_{H^1(\Omega_2)} \leq  M(\lambda) \varepsilon\sqrt{1+ \delta^2}
\end{eqnarray}
where $\displaystyle \delta = \frac{L_1}{L_1-L_0}$.
\end{theorem}
{\bf Proof}\\
The function $u|_{\Omega_2} - u_2$ is the solution of the system:
\begin{eqnarray*} 
\left\{
\begin{array}{lllllllll}
\displaystyle  - \Delta (u -u_2) = 0\quad \mbox{in}\quad \Omega_2 \\
\\
% v(x,z,0) = v_0(x,z)\\
% \\
\displaystyle \frac{\partial (u-u_2)}{\partial n}= 0 \quad \mbox{on}\quad \Gamma^2_{T}\\
\\
\displaystyle  \frac{\partial (u-u_2)}{\partial n} + \kappa (u-u_2)= 0 \quad \mbox{on}\quad \Gamma^2_{B}  \\
\\
 u-u_2 = 0 \quad \mbox{on} \quad \Gamma_R.
\end{array}
\right.
\end{eqnarray*}
By multiplying the first equation by $u-u_2$ and using the boundary conditions 
on $\Gamma_T^2 \cup \Gamma_B^2 \cup \Gamma_R$, we obtain:
\begin{equation} \label{eqn : etoile1}
\int_{\Omega_2}^{} |\nabla\left( u- u_2 \right)|^2 dx dz + \int_{\Gamma_B^2}^{} \kappa |u-u_2|^2 dx
 - \int_{\Gamma}^{} \frac{\partial (u-u_2)}{\partial n}\left( u - u_2 \right) dz  = 0. \tag{$*$}
\end{equation}
The integral term on $\Gamma$ is reformulated using the boundary condition 
$\displaystyle -\frac{ \partial u_2}{\partial x} + \lambda u_2 = -\frac{\partial u_1}{\partial x} + \lambda u_1$  satisfied by the limit $u_2$: 
\begin{eqnarray*}
\int_{\Gamma}^{} \frac{\partial (u-u_2)}{\partial n}\left( u - u_2 \right) dz  &=& - \int_{0}^{H}  \frac{\partial (u-u_2)}{\partial x}(L_0,z)\left( u - u_2 \right)(L_0,z) dz \\
 &= &  - \int_{0}^{H}  \frac{\partial u}{\partial x}(L_0,z)\left( u - u_2 \right)(L_0,z) dz \\
  & &   -\int_{0}^{H} \left(- \frac{\partial u_1}{\partial x}(L_0)  + \lambda u_1(L_0)  - \lambda u_2(L_0,z) \right)\left( u - u_2 \right) (L_0,z)dz \\
 &=&  \int_{0}^{H} \left(- \frac{\partial u}{\partial x}(L_0,z)  + \lambda u(L_0,z)   \right)\left( u - u_2 \right) (L_0,z)dz \\
     & &   -\int_{0}^{H} \left(- \frac{\partial u_1}{\partial x}(L_0)  + \lambda u_1(L_0)   \right)\left( u - u_2 \right)(L_0,z) dz \\
   && - \lambda\int_{0}^{H} \left(u - u_2 \right)^2(L_0,z) dz.
\end{eqnarray*}
$\bullet$ The first term reads:
\begin{eqnarray*}
\int_{0}^{H} \left(- \frac{\partial u}{\partial x}  + \lambda u   \right)(L_0,z)\left( u - u_2 \right) (L_0,z)dz  &=& 
\int_{0}^{H} \left(- \frac{\partial u}{\partial x}  + \lambda u\right) (L_0,z)\left( u(L_0,z) - \bar{u}(L_0) \right) dz \\
 && +  \int_{0}^{H} \left(- \frac{\partial u}{\partial x}  + \lambda u\right)(L_0,z)\left( \bar{u}(L_0) - u_1(L_0) \right) dz \\
 && + \int_{0}^{H} \left(- \frac{\partial u}{\partial x} + \lambda u \right)(L_0,z)\left( u_1(L_0) - u_2(L_0,z) \right) dz. \\
\end{eqnarray*}
Due to the relations (\ref{eqn : relation_ordre_2}) and (\ref{eqn : relation_ordre_2_bis}) and to the fact that $\displaystyle -\frac{\partial u}{\partial x}(L_0,z)  + \lambda u(L_0,z)= O(1) $, we deduce: 
$$
\int_{0}^{H} \left(- \frac{\partial u}{\partial x}  + \lambda u\right) (L_0,z)\left( u(L_0,z) - \bar{u}(L_0) \right) dz = O(\varepsilon^2)
$$
In the same way, if we assume that $L_0 < L_1$, so that 2-D effects are insignificant in $\Omega_1 \cap {\{ L_0 \leq x \leq L_1\}}$, and applying 
a similar asymptotic analysis as in the first section to the 2-D model defined in $\Omega_2$, we can deduce that:
\begin{eqnarray} \label{eqn : relation_u2_u1}
u_2(L_0,z) &=& \bar{u}_2(L_0) + O(\delta^2 \varepsilon^2) \nonumber \\
 &= & u_1(L_0) + O(\delta^2 \varepsilon^2),\quad \forall  z \in [0, H]
\end{eqnarray}
%{\bf Remark:}\\
So that:
$$
\int_{0}^{H} \left(- \frac{\partial u}{\partial x}  + \lambda u\right) (L_0,z)\left( u_1(L_0) - u_2(L_0,z) \right) dz = O(\delta^2 \varepsilon^2)
$$
Finally:
\begin{eqnarray} \label{eqn : dem_erreur1}
\int_{0}^{H} \left(- \frac{\partial u}{\partial x}  + \lambda u\right) (L_0,z)\left( \bar{u}(L_0) - u_1(L_0) \right) dz = 
H \left(- \frac{\partial \bar{u}}{\partial x}  + \lambda \bar{u}\right) (L_0)\left( \bar{u}(L_0) - u_1(L_0)\right)
\end{eqnarray}
%Ensuite le terme $\displaystyle -\int_{0}^{H} \left(- \frac{\partial u_1}{\partial x}(L_0)  + \lambda u_1(L_0)   \right)\left( u - u_2 \right) dz $ s'écrit :
$\bullet$ Since $u_1(L_0) = \overline{u}_2(L_0)$, the second term reads:
\begin{align} \label{eqn : etoile2}
-\int_{0}^{H} \left(- \frac{\partial u_1}{\partial x}(L_0)  + \lambda u_1(L_0)   \right)\left( u - u_2 \right)(L_0,z) dz &=
 -H \left(- \frac{\partial( u_1 - \bar u)}{\partial x}  + \lambda (u_1 -\bar{u}  )\right)(L_0)\left( \bar{u} - u_1 \right)(L_0)  \notag\\
&- H \left(- \frac{\partial \bar{u}}{\partial x}  
+ \lambda \bar{u}\right) (L_0)\left( \bar{u} - u_1\right)(L_0) \tag{$**$}
\end{align}
We reformulate the first term of the right. Note that the function $u_1 - \bar{u}$ satisfies the equation:
$$
-\frac{\partial^2 (u_1 - \bar{u})}{\partial x^2}(x) + a^2(u_1 -\bar{u})(x) = a^2(\bar{u}(x) - u(x,0)),\;\;\; \forall x \in (0, L_0)
$$
So that by multiplying this equation by $u_1 - \bar{u}$, after integration on $(0, L_0)$ and use of the boundary condition
$u_1(0) = \overline{u}(0)$, we obtain:
\begin{eqnarray*}
\int_{0}^{L_0} \left( \frac{\partial (u_1 - \bar{u})}{\partial x}\right)^2(x)dx  + a^2 \int_{0}^{L_0}\left( u_1 - \bar{u}\right)^2(x) dx 
- \frac{\partial (u_1 - \bar{u})}{\partial x}(L_0)(u_1 - \bar{u})(L_0)  =
\\
 \int_{0}^{L_0} a^2 (\bar{u}(x) - u(x,0))(u_1(x)-\bar{u})(x) dx
\end{eqnarray*}
thus:
\begin{eqnarray*}
- \frac{\partial (u_1 - \bar u)}{\partial x}(L_0)(u_1 - \bar{u})(L_0)  =
 \int_{0}^{L_0} a^2 (\bar{u}(x) - u(x,0))(u_1(x)-\bar{u}(x)) dx - \mathcal{A}_1(u_1-\bar{u}, u_1-\bar{u})
\end{eqnarray*}
where $\displaystyle\mathcal{A}_1(u_1-\bar{u}, u_1-\bar{u}) = \int_{0}^{L_0} \left( \frac{\partial (u_1 - \bar{u})}{\partial x}\right)^2(x)dx  + a^2 \int_{0}^{L_0}\left( u_1 - \bar{u}\right)^2(x) dx $.\\
And then (\ref{eqn : etoile2}) becomes:
\begin{eqnarray} \label{eqn : dem_erreur3}
-\int_{0}^{H} \left(- \frac{\partial u_1}{\partial x}(L_0)  + \lambda u_1(L_0)   \right)\left( u - u_2 \right) dz 
&=&  H\int_{0}^{L_0} a^2 (\bar{u}(x) - u(x,0))(u_1(x)-\bar{u}(x)) dx \nonumber\\
&&- H\mathcal{A}_1(u_1-\bar{u}, u_1-\bar{u})
+ \lambda H \left( \bar{u} - u_1 \right)^2(L_0) \nonumber \\
&& - H \left(- \frac{\partial \bar{u}}{\partial x}  
+ \lambda \bar{u}\right) (L_0)\left( \bar{u} - u_1\right)(L_0)
\end{eqnarray}
$\bullet$ To recap, the boundary term on $\Gamma$ in (\ref{eqn : etoile1}) becomes:
\begin{eqnarray*}
\int_{\Gamma}^{} \frac{\partial (u-u_2)}{\partial n}\left( u - u_2 \right) dz 
   &=& O(\varepsilon^2) + O(\delta^2 \varepsilon^2) + H\int_{0}^{L_0} a^2 (\bar{u}(x) - u(x,0))(u_1(x)-\bar{u}) dx \\
  && -H \mathcal{A}(u_1-\bar{u}, u_1-\bar{u})
   + \lambda H \left( \bar{u} - u_1 \right)^2(L_0)\\
   && - \lambda\int_{0}^{H} \left(u - u_2 \right)^2 dz 
\end{eqnarray*}
We first observe that: 
$$
\int_{0}^{L_0} a^2 (\bar{u}(x) - u(x,0))(u_1(x)-\bar{u}(x)) dx  \leq \frac{a^2}{2} \int_{0}^{L_0}  (\bar{u}(x) - u(x,0))^2dx + \frac{a^2}{2} \int_{0}^{L_0}  (u_1(x)-\bar{u}(x))^2dx.  
$$
It follows that:
\begin{eqnarray*}
\int_{\Gamma}^{} \frac{\partial (u-u_2)}{\partial n}\left( u - u_2 \right) dz & \leq &  C(1+ \delta^2 )\varepsilon^2) +
H \frac{a^2}{2} \int_{0}^{L_0}  (\bar{u}(x) - u(x,0))^2dx \\
&& - H\mathcal{A}_1(u_1-\bar{u}, u_1-\bar{u})
+H \frac{a^2}{2} \int_{0}^{L_0}  (u_1(x)-\bar{u}(x))^2dx \\
&& + \lambda H \left( \bar{u} - u_1 \right)^2(L_0) - \lambda\int_{0}^{H} \left(u - u_2 \right)^2(L_0,z) dz
\end{eqnarray*}
where $C$ is a positive constant. \\
Then we have:
$$\displaystyle - \mathcal{A}_1(u_1-\bar{u}, u_1-\bar{u})
+ \frac{a^2}{2} \int_{0}^{L_0}  (u_1(x)-\bar{u}(x))^2dx \leq 0$$
and finally, using the definition of $\overline{u}(L_0)$ and the relation $u_1(L_0) = \overline{u}_2(L_0)$, we obtain:
\begin{eqnarray*}
\lambda H \left( \bar{u} - u_1 \right)^2(L_0) - \lambda\int_{0}^{H} \left(u - u_2 \right)^2 dz &=& \lambda H \left( \frac{1}{H}\int_{0}^{H}(u -u_2)dz  \right)^2- \lambda\int_{0}^{H} \left(u - u_2 \right)^2 dz\\
 &= & \lambda \frac{1}{H}  \left( \int_{0}^{H}(u -u_2)dz  \right)^2 - \lambda\int_{0}^{H} \left(u - u_2 \right)^2 dz\\
 &\leq& \lambda \frac{1}{H}\left( \int_{0}^{H}1 dz  \right)\left( \int_{0}^{H}(u -u_2)^2dz  \right) - \lambda\int_{0}^{H} \left(u - u_2 \right)^2 dz\\
 &\leq& 0.
\end{eqnarray*}
We now come back to (\ref{eqn : etoile1}), which gives:
\begin{eqnarray*}
\int_{\Omega_2}^{} |\nabla\left( u- u_2 \right)|^2 dx dz + \int_{\Gamma_B^2}^{} \kappa |u-u_2|^2 dx 
& \leq & M(1 + \delta^2) \varepsilon^2
\end{eqnarray*}
%car $\displaystyle - \mathcal{A}_1(u_1-\bar{u}, u_1-\bar{u})
%+ \frac{a^2}{2} \int_{0}^{L_0}  (u_1(x)-\bar{u}(x))^2dx \leq 0$.\\
and thus:
$$
\int_{\Omega_2}^{} |\nabla\left( u- u_2 \right)|^2 dx dz \leq M (1 + \delta^2) \varepsilon^2
$$
Where $M$ denotes a positive constant.\\
Finally, due to the fact that $u - u_2 = 0$ on $\Gamma_R$, and by using Poincar\'e inequality we can deduce the inequality (\ref{eq : error_control}).$\Box$\\
\begin{remark}~\\
\begin{itemize}
\item This proposition fails if we choose the interface position in a  zone where 2-D effects are significant. In this case relation (\ref{eqn : relation_u2_u1})
is no more available.
\item The right term of (\ref{eq : error_control}) is also an upper bound of $\|u_2^{\lambda_1} - u_2^{\lambda_2}\|_{H^1(\Omega_2)}$
 for all $\lambda_1$ and $\lambda_2$ positive.
\end{itemize}
\end{remark}
\section{Numerical results}
The test cases presented in this section illustrate the coupling method of 1-D and 2-D elliptic equations  based on Schwarz algorithm. All the computations have been done using the software package Freefem++ \cite{Freefem}, with a $P_2$ finite element discretization. \\
In the first part of this section, the two test cases will be described in details. In the second and third parts, we will focus on one hand on the Scharwz algorithm convergence and on the other hand on the comparison of the coupled solution with the reference solution in order to enlight the theoretical results obtained in the previous paragraphs.

\subsection{Description of the test cases} 
\subsubsection{Test \#1:}
\begin{figure}[h!]
\begin{center}
\includegraphics[width=14cm]{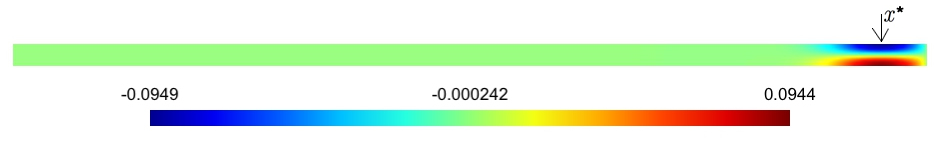}
 \caption{2-D reference solution for the first test case where $L=20$, $x^*=19$, $H=0.5$ and $\kappa=0.001$.}
 \label{fig: global solution test1 }
\end{center}
\end{figure}

The first test case is concerned with the solution of the 2-D problem (\ref{eq:full2-D}) where the domain is a rectangle $\Omega=[0, L] \times [0, H]$ which is assumed to be uniformly shallow: $H \ll L$. 
Let us consider that the right-hand side term $F(x,z)$ of the full 2-D problem is  $\displaystyle F(x,z) = m\exp(-(x-x^*)^2) \sin(\frac{2\pi z}{H})$, where $x^* < L$.\\
The global reference solution $u$ is displayed in Figure \ref{fig: global solution test1 }.

We notice that the 2-D effects are due to the particular form of the forcing term $F$ and are located around $x^*$.\\ 

\begin{figure}[h!]
\begin{center}
\subfloat[Computational domain for the 2-D reference model]
{\includegraphics{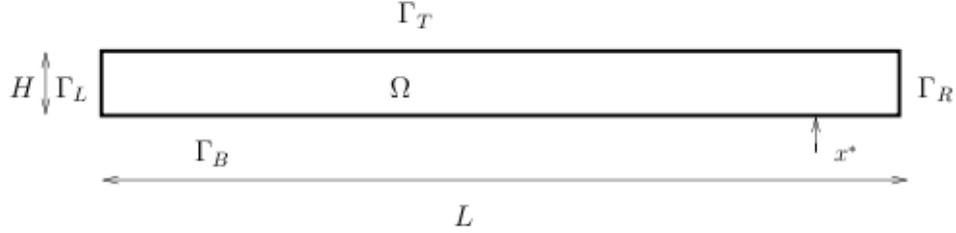} \label{fig : rectangulaire_entier}}\\
%{\input{Figs/domaine_rectangulaire.pdf_t} \label{fig : rectangulaire_entier}}\\
 \subfloat[Computational domain for the 1-D/2-D reduced model]
{\includegraphics{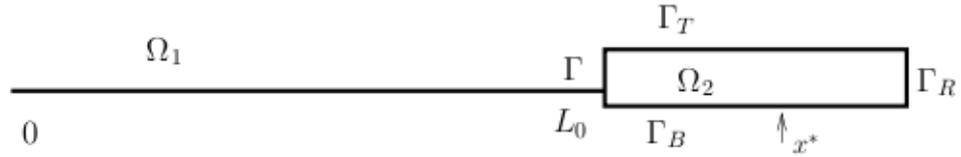} \label{fig : rectangulaire_couplage}} 
 \caption{Computational domains for both the reference and reduced models in test case \#1. For the reduced model (b), the  1-D/2-D interface $\Gamma$ is located in $x=L_0$.}
  \label{fig: domain test1 }
\end{center}
\end{figure}

Now let us define the coupled model. \\
The interface is located at $x=L_{0}< x^*$ as shown in Figure \ref{fig: domain test1 }. In the part of the domain $\Omega_{1}=[0, L_{0}] \times [0, H]$, we assume \textit{a priori} that the 2-D effects are negligible and consequently we replace the full 2-D equations by the 1-D model (see (\ref{eq:1-Dmodel})).  \\

\subsubsection{Test \#2:}
\begin{figure}[h!]
\begin{center}
\subfloat[Computational domain for the 2-D reference model]{
\includegraphics{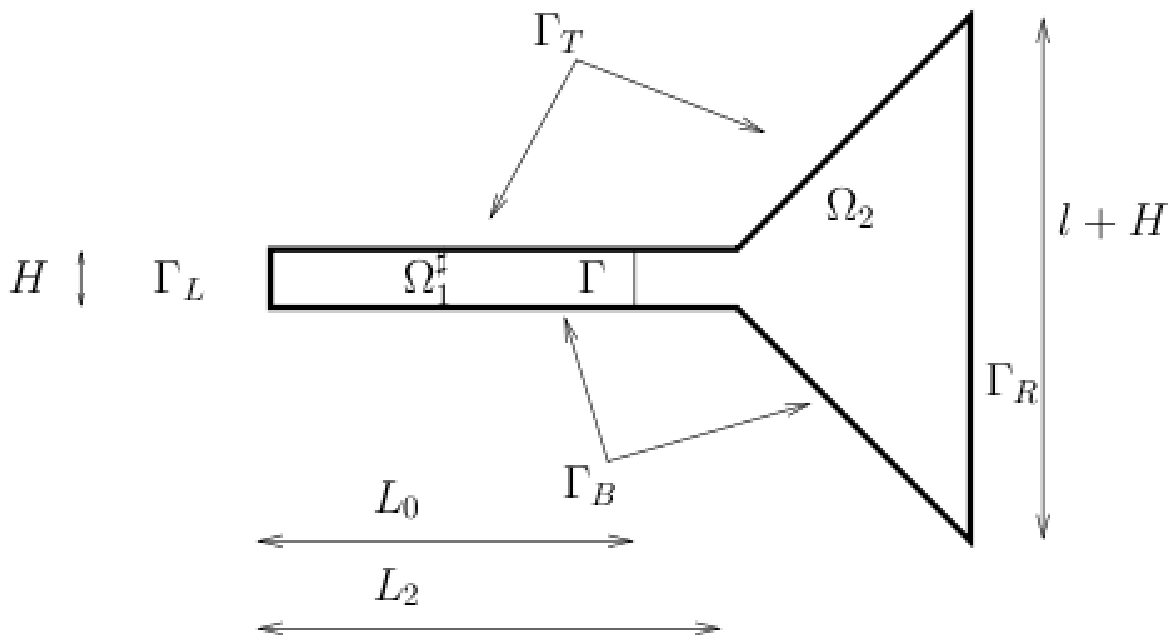}\label{fig : entonnoir}}\\
\subfloat[Computational domain for the 1-D/2-D reduced model]{
\includegraphics{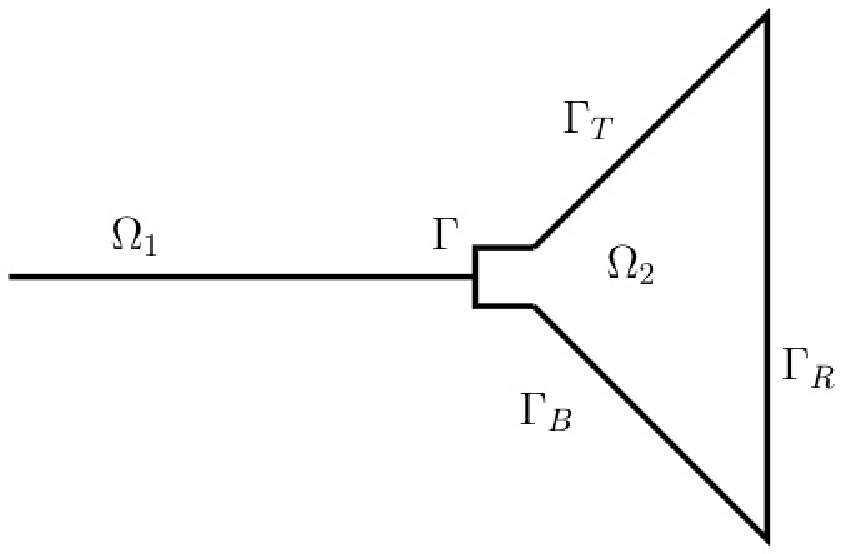}\label{fig : entonnoir_couplage}}
 \caption{Computational domains for both the reference and reduced models in test case \#2. For the reduced model (b), the  1-D/2-D interface $\Gamma$ is located in $x=L_0$.}
 \label{fig: domain test2 }
\end{center}
\end{figure}

In this second test case, the 2-D effects are due to the funnel-shaped geometry of the domain (see Figure \ref{fig : entonnoir}), and the forcing term is constant ($=1$).   \\

\begin{figure}[h!]
\begin{center}
\includegraphics[width=10cm]{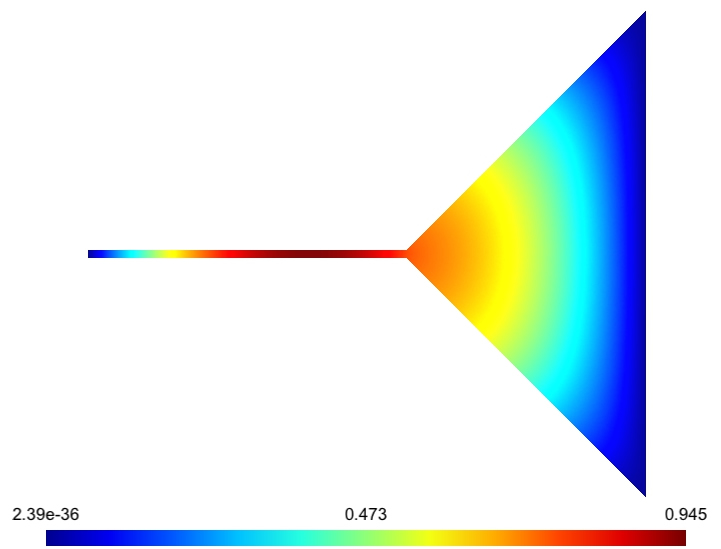}
 \caption{2-D reference solution for the second test case where $L=2$, $H=0.05$, $l=3$ and $\kappa=0.001$}
 \label{fig: global solution test2 }
\end{center}
\end{figure}
The reference solution in the whole domain is displayed in Figure \ref{fig: global solution test2 }.\\
The coupled model is defined by splitting the domain in two parts. The interface $\Gamma$ is located at $x = L_0$,  $0 < L_0 < L_2$ as shown in Figure \ref{fig: domain test2 }. 

\subsection{Convergence of the Schwarz algorithm}

\begin{figure}[h!]
\centering 
\subfloat{
\includegraphics[width=7cm]{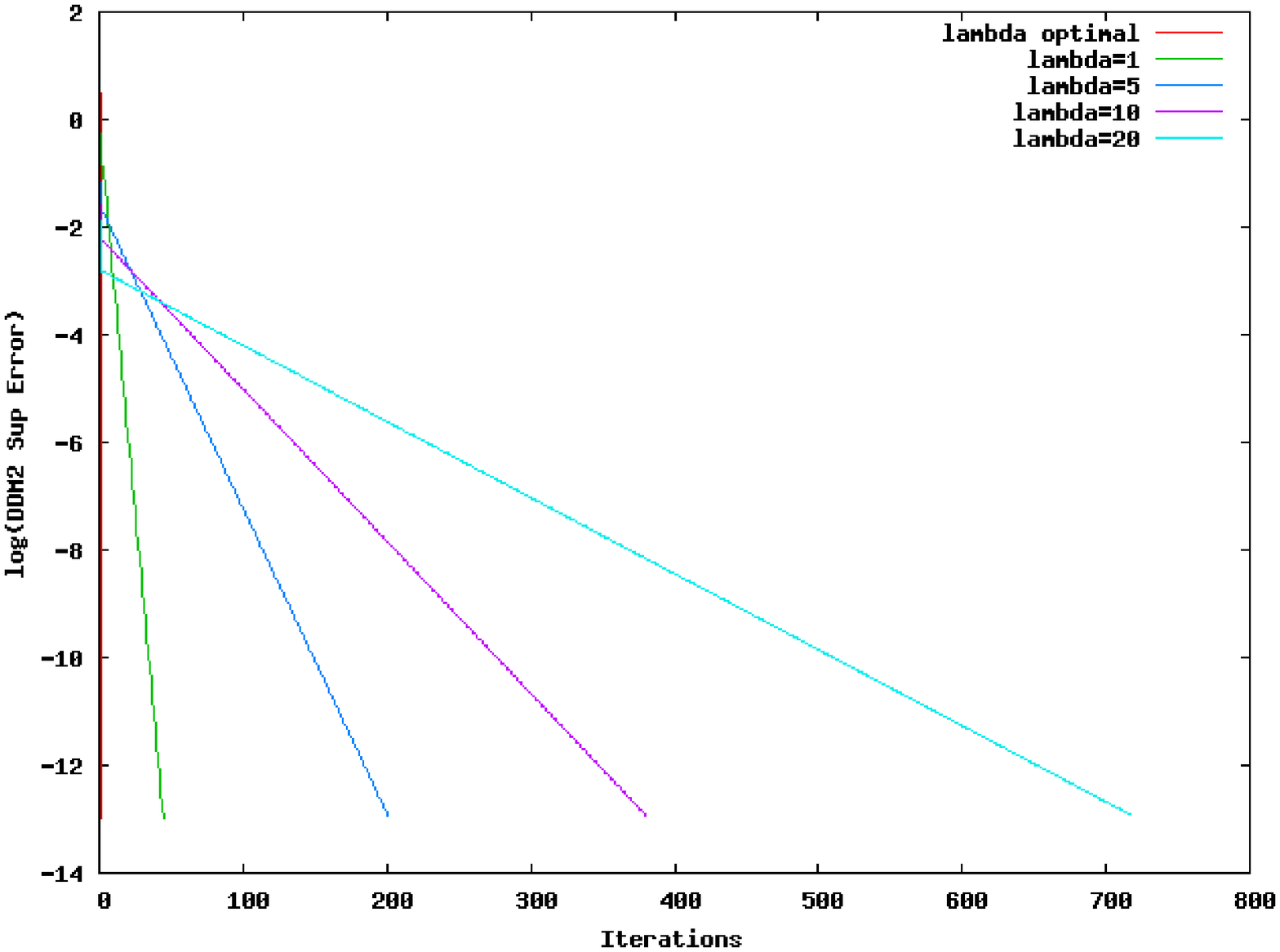} \label{fig: Schwarz_test_rect}
\includegraphics[width=7cm]{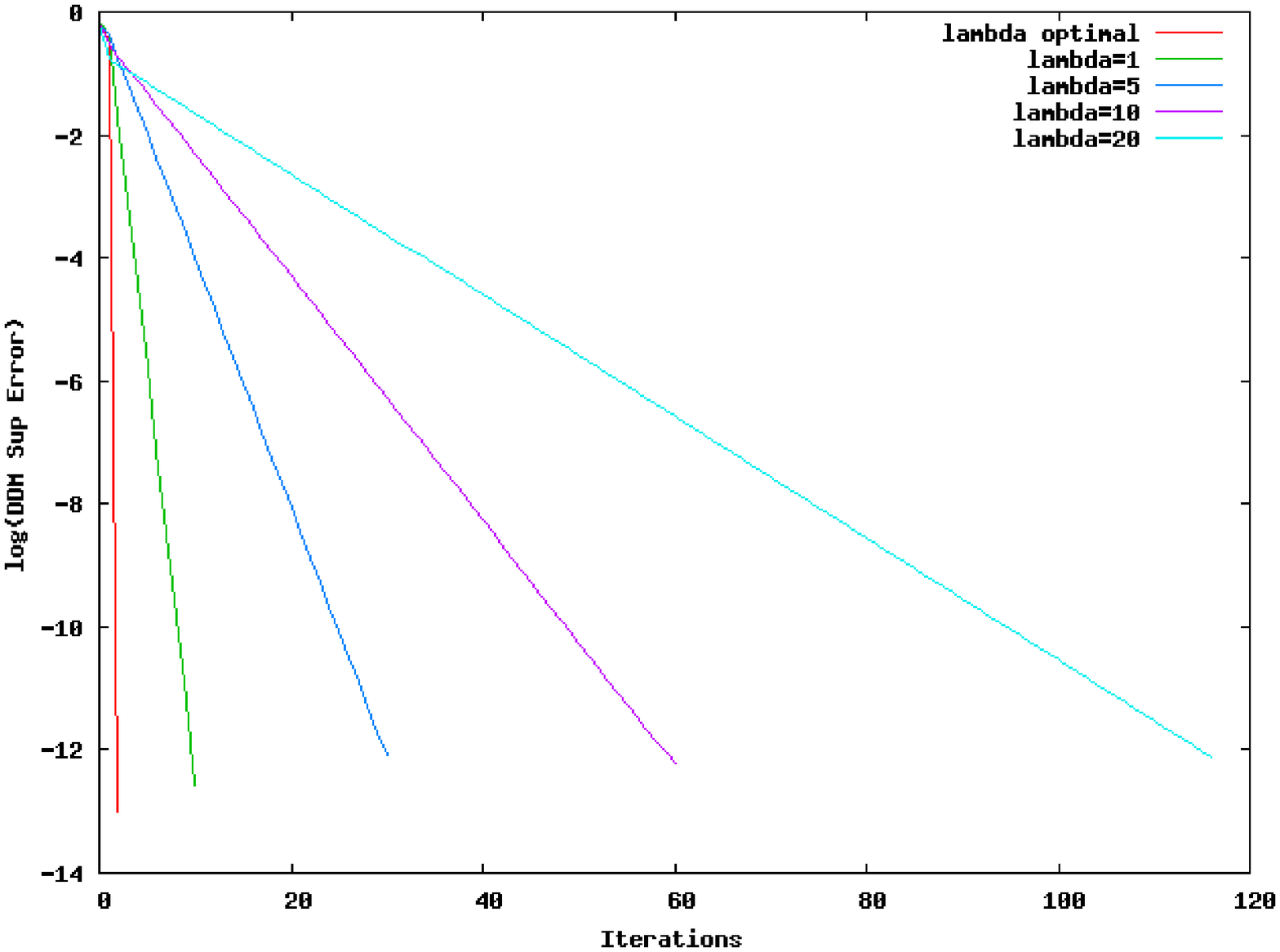}
} \label{fig: Schwarz_test_ent}
\caption{Convergence of Schwarz algorithm with various values of $\lambda$. Left: test case \#1 with $L=20$, $L_0=16$, $H=0.5$ and $\kappa=0.001$. Right: test case \#2 with $L=2$, $L_0=1.5$, $H=0.05$, $l=3$ and $\kappa=0.001$.}
\end{figure}

In this section we provide numerical results to assess the theoretical results of \S \ref{algorithm}. We are interested in illustrating the optimal convergence of Schwarz algorithm for the parameter $\lambda=\lambda_{opt}$. Figure \ref{fig: Schwarz_test_rect} shows the difference between the iterates of the Schwarz algorithm in $L^{\infty}$ norm for the two test cases.\\
As demonstrated in \S \ref{algorithm}, the Schwarz algorithm converges in two iterations for the optimal parametrer $\lambda=\lambda_{opt}$. It is important to notice that this result is independent of the interface location.\\

\subsection{Difference between the coupled solution and the full 2-D solution}
\begin{figure}[h!]
\centering 
\subfloat[\underline{Test case \# 1}: {\scriptsize  (left): $L=20$, $H=0.5$ and $\kappa=0.001$. (right): $L=20$, $L_0=14$ and $\kappa=0.01$.}]{
\includegraphics[width=7.5cm]{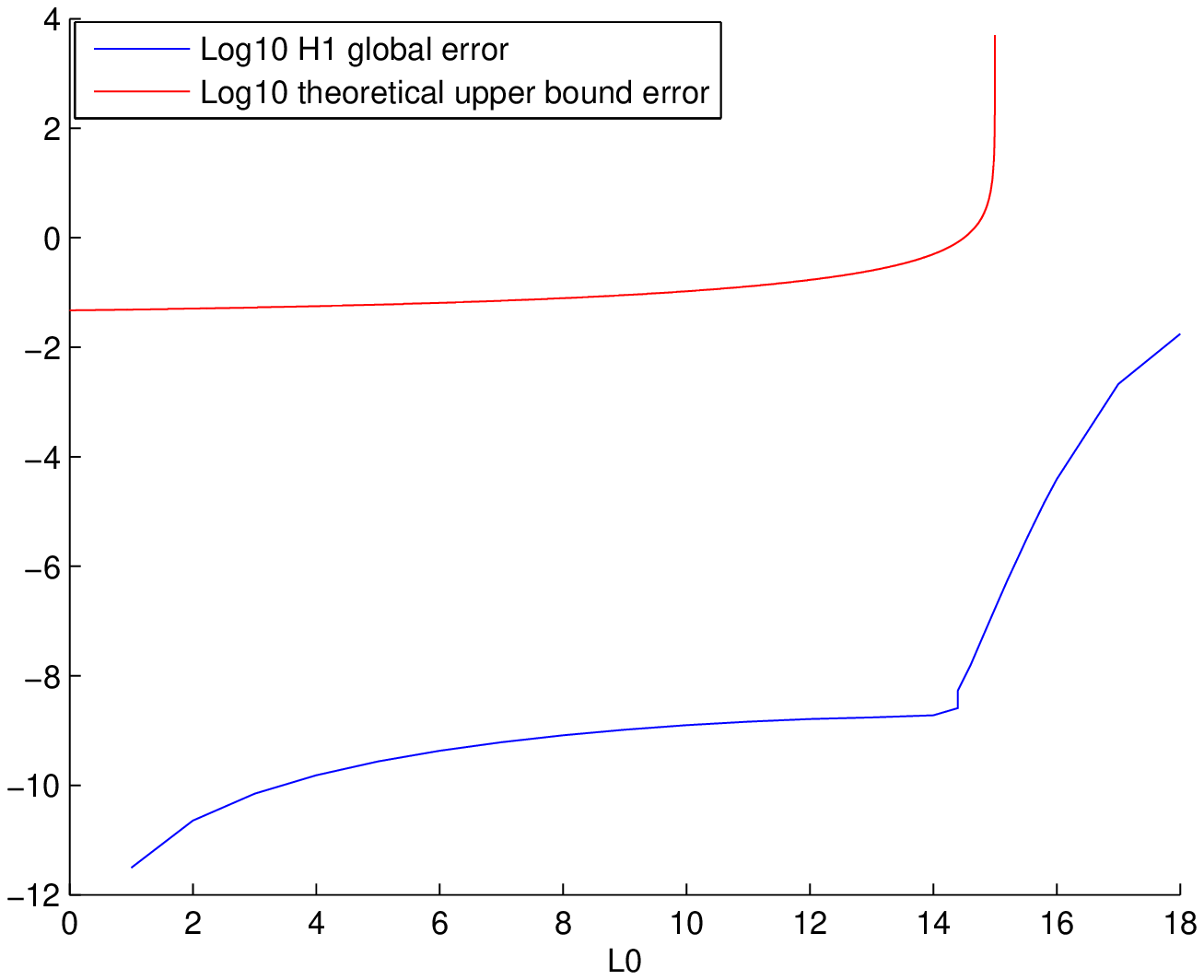} \label{fig: Erreur_sol_ref_rect1}
\includegraphics[width=7cm]{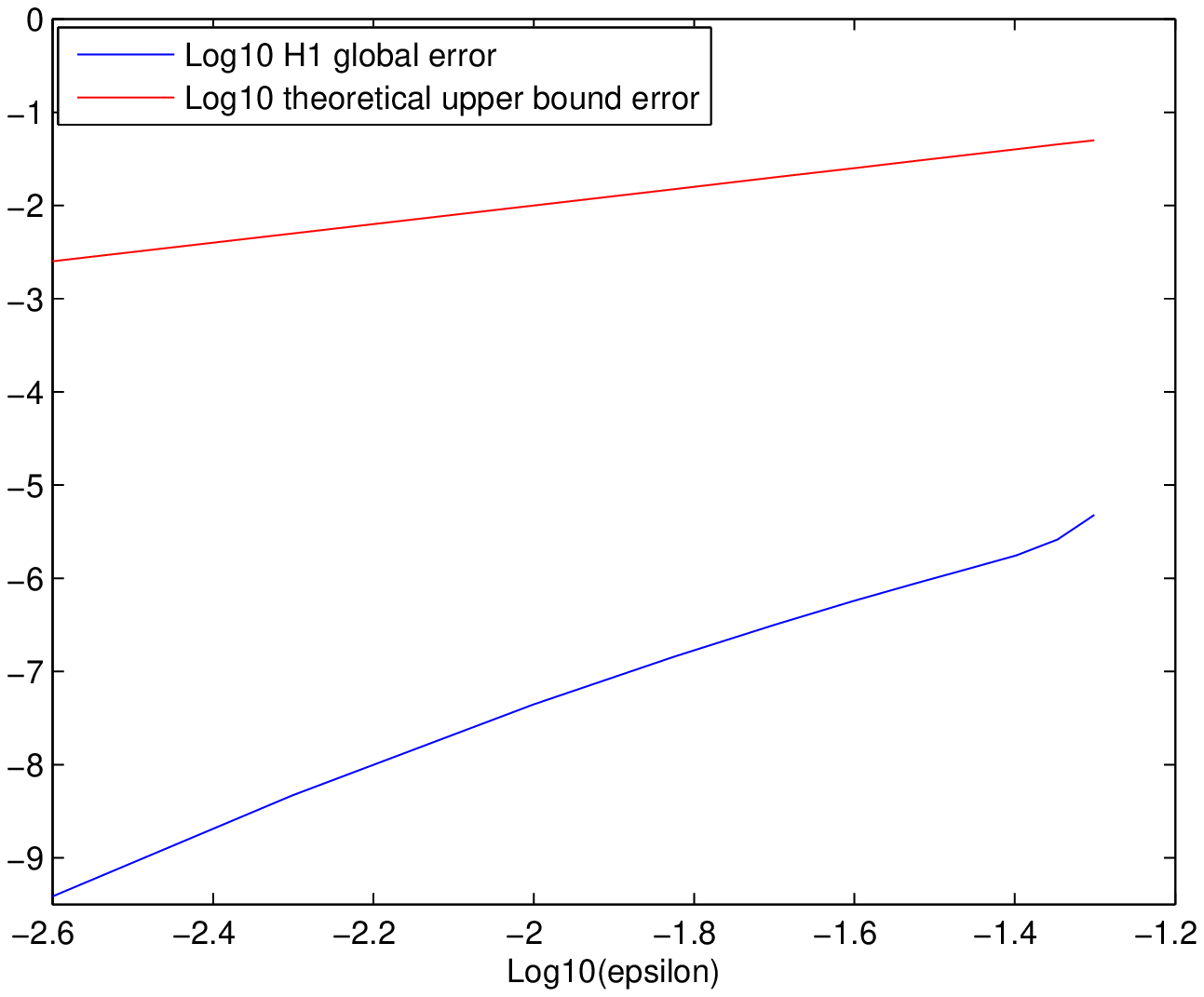} \label{fig: Erreur_sol_ref__eps_rect1}
}\\
\subfloat[\underline{Test case \# 2}: {\scriptsize  (left): $L=2$, $H=0.05$ and $\kappa=0.001$. (right): $L=2$, $L_0=1.5$ and $\kappa=0.001$.}]{
\includegraphics[width=7.5cm]{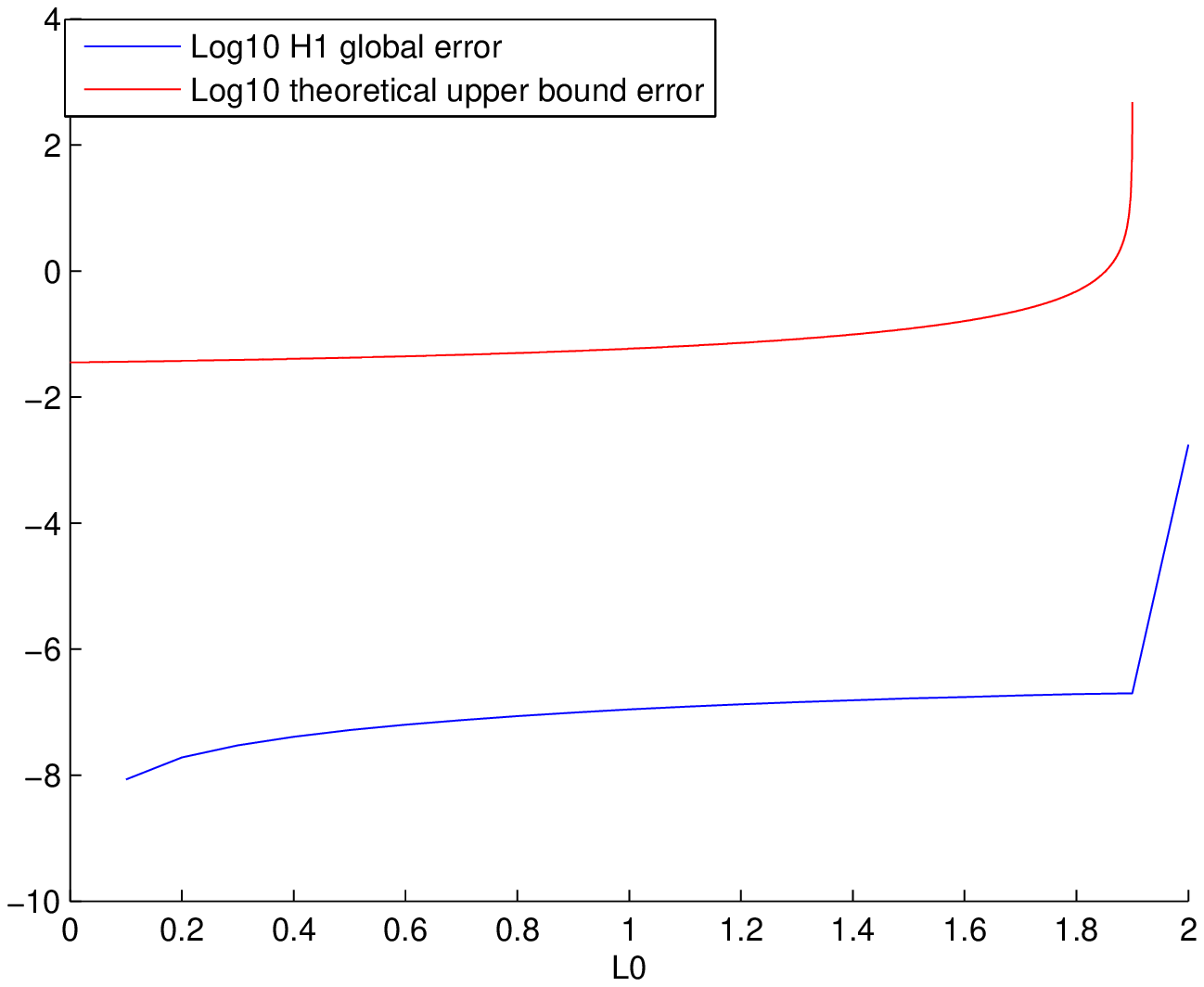} \label{fig: Erreur_sol_ref_ent1}
\includegraphics[width=7cm]{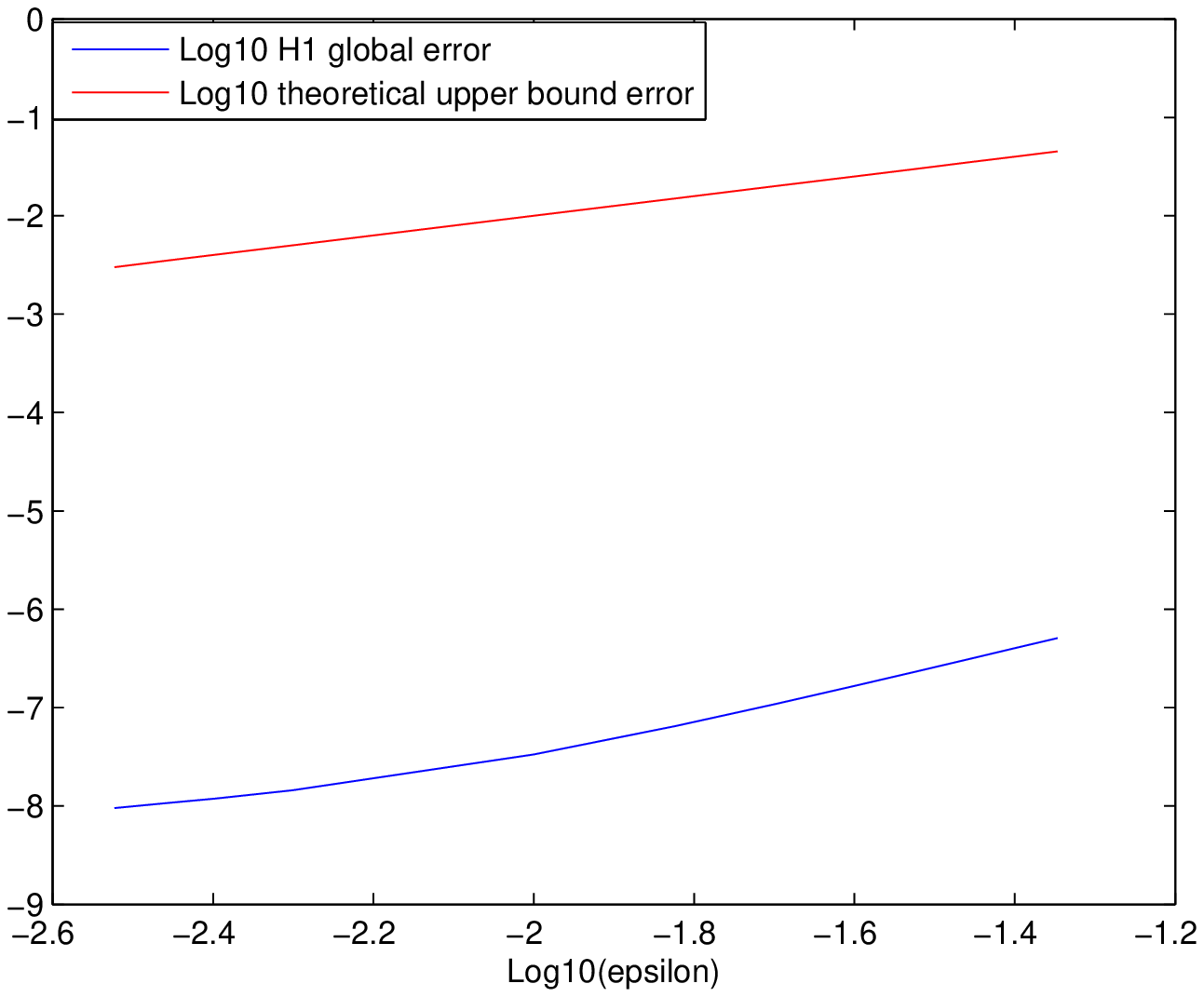} \label{fig: Erreur_sol_ref_eps_ent1}
}
 \caption{Relative error as a function of $L_0$ (left) and $\varepsilon$ (right) between the coupled solution and the 2-D reference solution in test case \#1 (top) and \#2 (bottom). 
 In the left column, the red curves correspond (for both test cases) to the RHS of estimate (\ref{eq : error_control}).}
 \label{fig:all}
\end{figure}

One important point in the analysis of the accuracy of the coupling procedure is the comparison of the coupled solution with the reference solution as a function of the interface location $x=L_{0}$. 
Contrarily to the classical domain decomposition problems, here there is a difference between the (converged) coupled solution and the full 2-D solution; this difference is due to the model reduction that is performed in the 1-D part of the domain. This difference depends on the location chosen to discriminate between 1-D and 2-D regions. Figures \ref{fig: Erreur_sol_ref_rect1} and \ref{fig: Erreur_sol_ref_ent1} left show the $H^1$ error between coupled and reference solutions as a function of the interface location for the two test cases. 
Figures \ref{fig: Erreur_sol_ref__eps_rect1} and \ref{fig: Erreur_sol_ref_eps_ent1}  right show the $H^1$ error in $\Omega_2$ between coupled and reference solutions as a function of $\displaystyle \varepsilon =\frac{H}{L}$ for the two test cases.\\

It is interesting to notice that for both test cases there is a discontinuity in the curve representing the error as a function of $L_{0}$ (see Figure \ref{fig:all}, left column). This discontinuity occurs both for the numerical difference between $u_2$ and $u_{\w_2}$ (black curve), and the theoretical curve (in red) corresponding to the right-hand-side of estimate (\ref{eq : error_control}). Indeed, if $L_{0}$ is greater than a certain threshold, the error grows very rapidly (and $\delta\longrightarrow\infty$ in estimate  (\ref{eq : error_control})). This could be an indication of the real (a priori unknown) value of $L_{1}$ (see discussion at the end of Section \ref{sec:derivation} above).

\newpage
\section{Conclusion}

In this paper we studied a linear boundary valued problem set in a 2-D domain, and assume that the solution may be approximated by a 1-D function in some part of the computational domain. We thus derive a reduced model that consists coupling a 1-D model (wherever we think it is legitimate) together with the original 2-D system (everywhere else).\\
The model reduction is performed thanks to a small aspect ratio hypothesis, with an integration in the shallow direction (we mimic the derivation of the shallow water equations). After this derivation we introduce an iterative method that couples the 1-D and 2-D systems and we prove some convergence results. One original aspect of this work 
is the particular attention that is paid to the location of the 1-D/2-D interface. These theoretical results are illustrated with numerical simulations that underline the importance of the interface position, but also the way 1-D and 2-D models are coupled (boundary conditions at this interface). All these aspects, that have been studied here with a linear model, will be considered in a forthcoming study of dimensionally heterogeneous modelling in fluid dynamics.

\section*{Acknowledgments}
This work was supported by the research department of the French national electricity company, 
\textit{EDF R}\&\textit{D}.

\tableofcontents

\end{document}